\documentclass[10pt]{article}                                                                                                                           
\usepackage{amssymb}                                                            
\usepackage{amsmath,amsthm} 
\usepackage{mathrsfs}                                                           
\usepackage{slashbox}  
\usepackage[all]{xy}
\usepackage[dvips]{graphicx}
\usepackage{wrapfig} 
\renewcommand{\t}{\tilde}
\newcommand{\ud}{\mathrm{d}}
\newcommand{\p}{\partial}
\newtheorem{prop}{Proposition}[section]
\newtheorem{col}{Corollary}[section]

\newtheorem{lemma}{Lemma}[section]
\newtheorem{theorem}{Theorem}[section]

\title{\bf{Quantum Group as \\
Semi-infinite Cohomology}}
\author{Igor B. Frenkel\footnote{igor.frenkel@yale.edu}\\
\vspace{5mm}Anton M. Zeitlin\footnote{anton.zeitlin@yale.edu, http://math.yale.edu/$\sim$az84 http://www.ipme.ru/zam.html}\\   
Department of Mathematics,\\
Yale University,\\
442 Dunham Lab, 10 Hillhouse Avenue,\\
New Haven, CT 06511}

\begin{document}
\maketitle
\begin{abstract}
We obtain the quantum group $SL_q(2)$ as semi-infinite cohomology of the Virasoro algebra with values in a tensor 
product of two braided vertex operator algebras with complementary central charges $c+\bar{c}=26$.
Each braided VOA is constructed from the free Fock space realization of the Virasoro algebra with 
an additional $q$-deformed harmonic oscillator degree of freedom. The braided VOA structure arises from the theory of
local systems over configuration spaces and it yields an associative algebra structure on 
the cohomology. We explicitly provide the four cohomology classes that serve as the generators of $SL_q(2)$
and verify their relations. We also discuss the possible extensions of our construction and its connection to 
the Liouville model and minimal string theory.
\end{abstract}
\tableofcontents
\section{Introduction}

Soon after the original discovery of the theory of quantum groups and their representations by Drinfeld \cite{D}
and Jimbo \cite{J}, several mathematicians and physicists have realized their profound relation to 
the representation theory of affine Lie algebras and conformal field theory \cite{mr}, \cite{fl}, \cite{gs}, \cite{rrr}, \cite{fw}, \cite{ffk}.
Eventually, this relation has been accomplished in the precise form of equivalence of certain tensor 
categories of representations \cite{KL}, \cite{F}, \cite{HL}. The nonstandard tensor product of the representations of
affine Lie algebras of the same level, motivated by two dimensional conformal field theory,
becomes natural in the context of vertex operator algebras (VOA) \cite{FHL}. For any simple Lie algebra 
$\mathfrak{g}$ let $\mathcal{C}_q$ be a category of type $I$ finite-dimensional representations of the quantum 
group $U_q(\mathfrak{g})$ and let $\mathcal{C}_k$ be a category of standard modules of the affine Lie algebra 
$\mathfrak{g}$. The simple objects of both categories are indexed by positive highest weights, which in the case 
$\mathfrak{g}=sl(2)$ can be identified with $\mathbb{Z}_+$. The equivalence of braided tensor categories $\mathcal{C}_q$ and $\mathcal{C}_k$, can be made transparent if one considers two other intermediate equivalent categories
\begin{equation}\label{eq:first}
\mathcal{C}_q\cong \mathcal{C}_\varkappa\cong \mathcal{C}_c\cong \mathcal{C}_k
\end{equation}
based, respectively, on the homology of configuration systems and certain representations of 
the $\mathcal{W}$-algebra, corresponding to $\mathfrak{g}$ (see \cite{dissert}).

The first isomorphism has been intensively studied by Varchenko et al (see \cite{varchenko},\cite{varchenko2}  and references therein);
the second isomorphism in the special case $\mathfrak{g}=sl(2)$, when $\mathcal{W}$-algebra is just the Virasoro algebra, is
implicit in the work of Feigin and Fuks \cite{FeFu}; the third isomorphism is a version of the quantum 
Drinfeld-Sokolov reduction developed in \cite{FeFr}.

The equivalence of representation categories points to a direct relation between the regular representations 
of the quantum group and affine Lie algebra, i.e. between Drinfeld's deformed algebra of functions on
the group and WZW conformal field theory. The latter does not have a vertex operator algebra structure that
would place it in the context of tensor categories, but it admits a remarkable modification 
that does have a VOA structure (\cite{frst} and references therein). Besides, it turns out that the central
charge of this modified regular VOA is precisely the one that yields a nonzero semi-infinite 
cohomology. A modified regular VOA can also be defined for $\mathcal{W}$-algebras at the critical central charge; 
in the special case of Virasoro algebra, the central charge is equal to 26, pointing to a relation with 
string theory \cite{frst}.

In fact, the relation between the theory of semi-infinite cohomology and string theory \cite{FGZ} has been
realized at about the time of the discovery of quantum groups and played an important role in the development
of both subjects. In particular, Lian and Zuckerman have shown in \cite{lz2} that 
the zero semi-infinite cohomology of vertex operator algebras inherits a natural structure of associative 
and commutative algebra, which was realized as a "ground ring'' in string theory \cite{wz}, \cite{wz2}.
In the case of modified regular VOA, the semi-infinite cohomology was identified as the center of 
the corresponding quantum group \cite{frst}.

In the present paper, we obtain the full quantum group $SL_q(2)$ via the semi-infinite cohomology.
Since the semi-infinite cohomology of any VOA is necessarily commutative we need to replace 
the modified regular VOA by a certain generalization, known as braided VOA \cite{frabs}, \cite{mr} (based on previous work \cite{tk}, \cite{ms}, \cite{ffk}). 
The latter is constructed from the tensor product of two braided VOA's with the complementary 
charges of precisely the same values as in the modified regular VOA, which ensures
nontriviality of the semi-infinite cohomology \cite{FGZ}.
In our paper, we treat in detail only the case of the Virasoro algebra, since it is 
the most interesting for pure mathematical reasons and because it has important applications in
physics. An extension of our construction to $\widehat{sl}(2)$ and other types of $\mathcal{W}$ and affine Lie 
algebras is straightforward though technically more difficult.
The key to our realization of a quantum group is an algebra isomorphism
\begin{equation}
H^{\frac{\infty}{2}+0}({\rm{Vir}},\mathbb{C}c,\mathbb{F}_c\otimes \mathbb{F}_{\bar{c}})\cong SL_q(2),
\end{equation}
where $c+\bar{c}=26$, and $q$ depends on $c$.
The braided VOA $\mathbb{F}_c$ (and similarly $\mathbb{F}_{\bar{c}}$) can be realized on the space
\begin{equation}\label{space}
\oplus_{\lambda\ge 0}(V_{\Delta(\lambda),c}\otimes V_\lambda),
\end{equation}
where $V_\lambda$ and $V_{\Delta(\lambda),c}$ are the corresponding simple modules in 
the equivalent braided tensor categories $\mathcal{C}_q$ and $\mathcal{C}_c$ of quantum algebra $U_q(sl(2))$ and the Virasoro algebra. The most transparent way 
to describe the braided VOA structure arises from the equivalence of both categories to 
the intermediate one $\mathcal{C}_\varkappa$ in (\ref{eq:first}). 
As we mentioned before, the category $\mathcal{C}_\varkappa$ is constructed from the homology of configuration spaces and it provides a geometric realization 
of the purely algebraic structure of $\mathbb{F}_c$. However, it is more convenient to consider a generic version $\tilde{\mathcal{C}}_\varkappa$ of $\mathcal{C}_\varkappa$, which also exists for the other braided categories in (\ref{eq:first}), and we also obtain their equivalences
\begin{equation}\label{cat2}
\tilde{\mathcal{C}}_q\cong \tilde{\mathcal{C}}_\varkappa\cong \tilde{\mathcal{C}}_c\cong \tilde{\mathcal{C}}_k.
\end{equation}
In these categories we allow the weights to be arbitrary complex numbers , therefore the simple objects in each category are indexed not by $\mathbb{Z}_+$ as in (\ref{eq:first}) but by $\mathbb{C}$. They are again simple highest weight modules with fixed central extensions  as in (\ref{eq:first}) though for the generic weights in $\mathbb{C}\backslash \mathbb{Z}_+$ they coincide with the Verma modules and the contragradient Verma modules. It is also convenient to realize them as certain Fock spaces. At the generic highest weights we have semisimplicity in all four categories of (\ref{cat2}), which allow us to define the structure coefficients of our generic categories in a complete parallel with the classical categories in (\ref{eq:first}). We verify that the structure coefficients of 
$\tilde{\mathcal{C}}_q$, $\tilde{\mathcal{C}}_\varkappa$, $\tilde{\mathcal{C}}_c$ at the generic weights are the same (the category $\tilde{\mathcal{C}}_k$ is not considered in this paper but the equivalence with other three categories is straightforward). Then, using the Fock space realization of the Verma modules we analytically continue certain 
structure coefficients that allow us to relate the structure coefficients in (\ref{cat2}) with the ones in (\ref{eq:first}). 
Note, that for $\mathfrak{g}=sl(2)$ the relation between (\ref{eq:first}) and (\ref{cat2}) can be restated in terms of analytic continuation of quantum $6j$-symbols and the corresponding identities. 
The full analysis of the categories in (\ref{cat2}) and their equivalences at the integral points is an interesting problem, the solution of which, however, is not necessary for the main goal of the present paper. 

The structure of the homology of configuration spaces that determines $\tilde{\mathcal{C}}_\varkappa$ and the isomorphism with the representation category of quantum group $\tilde{\mathcal{C}}_q$ was extensively studied by Varchenko in \cite{varchenko}, \cite{varchenko2} for an arbitrary type of Lie algebra. On the other hand, in the Fock space realization of the representations of Virasoro and $\hat{sl}(2)$ Lie algebras, the intertwining operators are given by integrals of vertex operators via certain cycles that precisely belong to $Hom$'s of the category $\tilde{\mathcal{C}}_\varkappa$, which yields the other two isomorphisms of (\ref{cat2}) from their generic counterparts by a limiting procedure that we explain in our paper. 
In particular, this leads to a realization of the space (\ref{space}) as a subspace of the Fock space
\begin{equation}
S(a_{-1},a_{-2},\dots)\otimes S(\beta)\otimes\mathbb{C}[\mathbb{Z}]
\end{equation}
with the natural action of the Heisenberg algebra
\begin{equation}
[a_m,a_n]=2\varkappa m\ \delta_{m+n,0}
\end{equation}
and the 1-dimensional $q$-deformed harmonic oscillator
\begin{equation}
\gamma\beta-q\beta\gamma=q^{-N}, \quad [N,\beta]=\beta, \quad [N,\gamma]=-\gamma.
\end{equation}
The verification of the axioms of braided VOA for $\mathbb{F}_c$ is essentially identical to the 
proof of the equivalence of these tensor categories (see \cite{dissert}).
Though the latter equivalence has been studied in many sources (though often in disguised form 
\cite{mr}, \cite{fl}), to make the paper self contained and explicit, we provide all the necessary details 
that are needed for the verification of the braided VOA structure. We also explicitly identify
the representatives of the semi-infinite cohomology classes that yield the generators $A,B,C,D$ of 
the quantum group $SL_q(2)$ and verify the defining relations
\begin{eqnarray}
&&AB=q^{-1}BA,\quad AC=q^{-1}CA,\quad BD=q^{-1}DB,\nonumber\\
&&CD=q^{-1}DC,\quad AD-DA=(q^{-1}-q)BC,\quad BC=CB,\nonumber\\
&&AD-q^{-1}BC=1.
\end{eqnarray}
Our results open new perspectives of the relation of semi-infinite cohomology and string theory. 
In fact, the coupling of $\mathbb{F}_c$ with $0\leqslant c\leqslant 1$ and $\mathbb{F}_{\bar{c}}$
with $25\leqslant c\leqslant 26$ can be interpreted as a coupling of a minimal model and the Liouville
model  with complementary central charges $c+\bar{c}=26$, which arises in the so-called minimal
string theory (see \cite{minimal} and references therein). This indicates that our realization of the quantum group as 
the semi-infinite cohomology might admit a geometric interpretation in terms of string theory providing 
a new link between two subjects. Thus, our construction might be viewed as a step towards 
an invariant geometric description of the untamed noncommutative structure of quantum group.

The paper is organized as follows. Section 2 is devoted to some basic facts about $U_q(sl(2))$. We remind statements useful in the following: the statements 
about Verma and dual Verma modules, and the relations between associated intertwiners. We also derive the 
polynomial ($q$-oscillator) realization for the intertwiner between dual Verma modules. The third section is devoted to the data we will work with throughout this paper. Namely, we consider the lattice of Fock modules and the braided VOA on this space associated with Feigin-Fuks realization of Virasoro algebra. In this section we remind basic facts about irreducible Virasoro modules 
for the generic values of central charge: we partly use the tools, which were introduced by Felder in the more complicated case of rational conformal field theory \cite{felder}. 

In Section 4, we study the geometry of local systems associated with the multivalued function corresponding to a certain correlator from the braided VOA constructed in Section 3. Our constructions are motivated by the heuristic constructions of Gomez and Sierra 
\cite{gs} and rigorous results of \cite{varchenko}, \cite{fkv}. These geometric considerations allow us to construct in Section 5 the braided vertex algebra of intertwiners between Fock spaces and then braided VOA on the space $\mathbb{F}_c$ (see (\ref{space})). 

In Section 6, we consider a certain "double" of the braided VOA from Section 5. Namely, we examine the structure of the braided VOA on the 
space $\mathbb{F}=\mathbb{F}_{c}\otimes \mathbb{F}_{\bar{c}}$. It appears that there is a hereditary ring structure on 
the semi-infinite cohomology of $H^{\frac{\infty}{2}+\cdot}(Vir, \mathbb{C}\mathbf{c},\mathbb{F})$. As we already mentioned above, one can 
explicitly calculate this semi-infinite cohomology on the zero level and show that it (as a space) 
coincides with $SL_q(2)$. Lian and Zuckerman introduced an associative product on the space of the semi-infinite cohomology of VOA. Applying the certain modification of the Lian-Zuckerman construction to our braided VOA, we reproduce the multiplicative structure of $SL_q(2)$ on the zero level of the semi-infinite cohomology space. 

In the last Section, we outline possible extensions of the results in this paper.

\section{$\mathbf{U_q(sl(2))}$, its representations and intertwining operators}
\setcounter{equation}{0}

{\bf {2.1. Basic facts and notations.}} Let $U_q(sl(2))$ be the Hopf algebra over $\mathbb{C}(q)$ with generators $E, F, q^{\pm H}$
and commutation relations:

\begin{eqnarray}
q^{\pm H}E&=&q^{\pm 2}Eq^{\pm H},\nonumber\\
q^{\pm H}F&=&q^{\mp 2}Fq^{\pm H},\nonumber \\
{[}E,F] &=&\frac{q^H-q^{-H}}{q-q^{-1}}.
\end{eqnarray}
The comultiplication is given by
\begin{eqnarray}
&&\Delta(q^{\pm H})=q^{\pm H}\otimes q^{\pm H}\nonumber\\
&&\Delta(E)=E\otimes q^H+1\otimes E,\nonumber\\
&&\Delta(F)=F\otimes 1 + q^{-H}\otimes F,
\end{eqnarray}
The universal R-matrix for $U_q(sl(2))$, which is an element of a certain completion of 
$U_q(sl(2))\otimes U_q(sl(2))$, is given by:
\begin{eqnarray}
R&=&C\Theta, \qquad C=q^{\frac{H\otimes H}{2}},\nonumber\\
\Theta&=&\sum_{k\geqslant 0}q{^{k(k-1)/2}\frac{(q-q^{-1})^k}{[k]!}E^k\otimes F^k},
\end{eqnarray}
where $[n]=\frac{q^n-q^{-n}}{q-q^{-1}}$ and $[n]!=[1][2]\ldots[n]$.

For any given pair $V,W$ of representations, R-matrix gives the following commutativity isomorphism:
$\check{R}=PR: V\otimes W \rightarrow W\otimes V$, where $P$ is a permutation: $P(v\otimes w)=w\otimes v$.

We denote by $M_\lambda$ the Verma module with highest weight $\lambda\in \mathbb{C}$. We will say that the weight $\lambda$ is $generic$, if $\lambda\notin \mathbb{Z}$. 
In the case $\lambda\in\mathbb{Z_+}$ one obtains an irreducible finite dimensional representation 
$V_\lambda$ (of dimension $\lambda+1$) by means of quotient:
\begin{equation}\label{factverma}
V_\lambda=M_\lambda/\langle F^{\lambda+1}v_\lambda\rangle,
\end{equation}
where $v_\lambda$ is the vector corresponding to the highest weight in $M_\lambda$.

Let us define an algebra anti-automorphism 
$\tau:U_q(sl(2))\rightarrow U_q(sl(2))$
by means of
$\tau(E)=F{q^H}, \tau(F)=Eq^{-H}, \tau(q^H)=q^H, \tau(ab)=\tau(b)\tau(a)$.
Then $\tau$ is a coalgebra automorphism:
$(\tau\otimes\tau)\Delta(x)=\Delta(\tau(x))$ and $\tau(R)=R^{21}=P(R)$,
where $P(a\otimes b)=b\otimes a$.
For every module M, let the contragradient module
$M^c$ be the restricted dual to M with the action of $U_q(sl(2))$ given by
\begin{equation}
<gv^*,v>=<v^*,\tau(g)v>,\quad v\in M, v^*\in M^c, g\in U_q(sl(2)).
\end{equation}
Note, that $(M_1\otimes M_2)^c\cong M_1^c\otimes M_2^c$ and for $\lambda\in\mathbb{Z}_+$
we have $V_\lambda^c\cong V_\lambda$. From (\ref{factverma}) we have an embedding for $\lambda\in \mathbb{Z}_+$:
\begin{equation}
V_\lambda\subset M_\lambda^c,
\end{equation}
which leads to the following exact sequence:
\begin{equation}\label{ex}
0\rightarrow V_\lambda\rightarrow M_\lambda^c\rightarrow V_{-\lambda-2}\rightarrow 0,
\end{equation}
and for $\lambda\in\mathbb{Z}_{\leqslant -1}$ we have $M_\lambda^c\cong V_\lambda$.\\

\noindent{\bf {2.2. Intertwining operators}}. In the following the intertwiners for Verma modules and their contragradient  counterparts will play a crucial role. 
Namely, in this paper we will consider the elements of 
$Hom(M_\mu^c\otimes M_\lambda^c, M_\nu^c)$ and its dual $Hom(M_\nu, M_\mu\otimes M_\lambda)$, such that $\lambda,\mu,\nu \in \mathbb{C}$ and $\mu+\lambda-\nu\in 2\mathbb{Z}$. In the generic case if $\nu\le \mu+\lambda$ then $dim Hom(M_\nu, M_\mu\otimes M_\lambda)=1$ and otherwise 
$dim Hom(M_\nu, M_\mu\otimes M_\lambda)=0$. 
We will use the following notation for the intertwining operators from 
$Hom(M_\mu^c\otimes M_\lambda^c, M_\nu^c)$:
\begin{equation}
\Phi^\nu_{\mu\lambda}(\cdot\otimes\cdot): M_\mu^c\otimes M_\lambda^c\rightarrow M_\nu^c
\end{equation}
and the ones from $Hom(M_\nu, M_\mu\otimes M_\lambda)$:
\begin{equation}
\Phi^{\mu\lambda}_\nu(\cdot): M_\nu\rightarrow M_\mu\otimes M_\lambda.
\end{equation}
 We will also need intertwiners from $Hom(V_\nu, V_\mu\otimes V_\lambda)$,  where 
$\mu,\nu,\lambda\in \mathbb{Z}_+$ and their dual from $Hom(V_\mu\otimes V_\lambda, V_{\nu})$. They can be reconstructed from the intertwiners above by means of the following projections/embeddings:
\begin{eqnarray}
&& M_\nu\xrightarrow{\Phi^{\mu\lambda}_\nu} M_\mu\otimes M_\lambda\xrightarrow{P_{\mu}\otimes P_{\lambda}}V_\mu\otimes V_\lambda,\nonumber\\
&& V_\mu\otimes V_\lambda\xrightarrow{i_{\mu}\otimes i_{\lambda}}M_\mu^c\otimes M_\lambda^c\xrightarrow{ \Phi^\nu_{\mu\lambda}}M_\nu^c,
\end{eqnarray}
where $P_{\xi}$ is a standard projection on the irreducible module from the corresponding Verma module, and $i_{\xi}$ is an embedding of the finite dimensional irreducible module into contragradient Verma module. It is clear that the first expression gives the element  from $Hom(V_\nu, V_\lambda\otimes V_\mu)$ and the second one corresponds to $Hom(V_\mu\otimes V_\lambda, V_\nu)$. 
Similarly, one can construct intertwiners from $Hom(M_\mu\otimes V_\lambda, M_\nu)$ and 
$Hom(M_\mu\otimes V_\lambda, M_\nu)$. 
Let us denote the elements of $Hom(V_{\nu},V_{\mu}\otimes V_{\lambda})$ and $Hom(V_{\mu}\otimes V_{\lambda},V_{\nu})$ as $\phi_{\nu}^{\mu\lambda}$ and 
$\phi^{\nu}_{\mu\lambda}$ correspondingly. 
It is known that there exist identifications
\begin{eqnarray}\label{homsing}
 &&Hom(M_{\nu},M_{\mu}\otimes M_{\lambda})\cong Sing_{\nu}(M_\mu\otimes M_{\lambda}),\\
&& Hom(V_{\nu}, V_{\mu}\otimes V_{\lambda})\cong Sing_{\nu}(V_{\mu}\otimes V_{\lambda}),
\end{eqnarray}
where $Sing_{\nu}$ denotes the space of singular vectors of the weight $\nu$. The explicit form of the above isomorphism is given by the following map:
\begin{equation}
\Phi^{\mu\lambda}_\nu\to \Phi^{\mu\lambda}_\nu (v_{\nu}),
\end{equation}
where $v_{\nu}$ is the highest weight vector in $M_{\nu}$. 

In the case of the generic values of the weights of appropriate modules, we have the following 
Proposition, expressing the bilinear relations for the intertwining operators.
 
\begin{prop}\label{intertwiner0}
Let $\lambda_i$ $(i=0,1,2,3)$ be generic. 
Then there exists an invertible operator 
\begin{displaymath}
B^M\left[\begin{array}{cc}
\lambda_0 \ \lambda_1\\
\lambda_2 \ \lambda_3
\end{array}
\right]
\end{displaymath}
such that the following diagram is commutative:

\begin{displaymath}
\xymatrix{
{\begin{array}{c}
\oplus_\rho \big(Hom(M_\rho,M_{\lambda_1}\otimes M_{\lambda_2})\\
\otimes Hom(M_{\lambda_0}, M_\rho\otimes  M_{\lambda_3})\big)
\end{array}
\ar[d]_{i}}
\ar[r]^{\qquad\qquad
B^M\left[\begin{array}{cc}
\lambda_0 \lambda_1\\
\lambda_2 \lambda_3
\end{array}
\right]
\qquad\qquad}&  
{\begin{array}{c}
\oplus_\xi \big(Hom(M_\xi,M_{\lambda_1}\otimes M_{\lambda_3})\\
\otimes Hom(M_{\lambda_0}, M_{\xi}\otimes M_{\lambda_2})\big)
\end{array}
\ar[d]_{i}
}
\\
Hom(M_{\lambda_0}, M_{\lambda_1}\otimes M_{\lambda_2}\otimes M_{\lambda_3})
\ar[r]^{PR}& 
Hom(M_{\lambda_0}, M_{\lambda_1}\otimes M_{\lambda_3}\otimes M_{\lambda_2}) ,
}
\end{displaymath}
where  $\rho\in \{\lambda_1+\lambda_2-2k$, $k\in\mathbb{Z}_+\}$,  
$\xi\in \{\lambda_1+\lambda_3-2k$, $k\in\mathbb{Z}_+\}$, and $i$ is an isomorphism.
\end{prop}
The proof follows from the complete reducibility of the tensor product Verma modules in the case of generic highest weights. 

Using the notation we introduced above, one can write the statements of Proposition \ref{intertwiner0} as follows:
\begin{eqnarray}\label{int1}
(1\otimes PR)\Phi_\rho^{\lambda_1\lambda_2}\Phi_{\lambda_0}^{\rho\lambda_3}=
\sum_{\xi}B^M_{\rho\xi}\left[\begin{array}{cc}
\lambda_0 & \lambda_1\\
\lambda_2 & \lambda_3
\end{array}
\right]
\Phi_\xi^{\lambda_1\lambda_3}\Phi_{\lambda_0}^{\xi\lambda_2},
\end{eqnarray}
where $B^M_{\rho\xi}\left[\begin{array}{cc}
\lambda_0 & \lambda_1\\
\lambda_2 & \lambda_3
\end{array}
\right]$ are the matrix coefficients of the operator $B^M$, which depend on the normalization of the intertwining operators $\Phi_\nu^{\lambda\mu}$. We will fix such normalization later. 

For the dual intertwiners $\Phi_{\mu\nu}^\lambda$,  we have similar identity:
\begin{eqnarray}\label{int2}
\Phi_{\rho\lambda_3}^{\lambda_0}\Phi_{\lambda_1\lambda_2}^\rho(1\otimes PR)=
\sum_{\xi}B^M_{\xi\rho}\left[\begin{array}{cc}
\lambda_0 & \lambda_1\\
\lambda_2 & \lambda_3
\end{array}
\right]
\Phi_{\xi\lambda_2}^{\lambda_0}\Phi_{\lambda_1\lambda_3}^\xi.
\end{eqnarray}
In the case of integer weights, the following Proposition gives the bilinear algebraic relations between the compositions of the intertwiners of finite-dimensional modules (see e.g. \cite{efk}).

\begin{prop}\label{intertwiner1}
Let $\lambda_i\in \mathbb{Z}_+$ $(i=0,1,2,3)$ . 
Then there exists an invertible operator 
\begin{displaymath}
B^V\left[\begin{array}{cc}
\lambda_0 \ \lambda_1\\
\lambda_2 \ \lambda_3
\end{array}
\right]
\end{displaymath}
such that the following diagram is commutative:

\begin{displaymath}
\xymatrix{
{\begin{array}{c}
\oplus_\rho \big(Hom(V_\rho,V_{\lambda_1}\otimes V_{\lambda_2})\\
\otimes Hom(V_{\lambda_0}, V_\rho\otimes  V_{\lambda_3})\big)
\end{array}
\ar[d]_{i}}
\ar[r]^{\quad
B^V\left[\begin{array}{cc}
\lambda_0 \lambda_1\\
\lambda_2 \lambda_3
\end{array}
\right]
\qquad}&  
{\begin{array}{c}
\oplus_\xi \big(Hom(V_\xi,V_{\lambda_1}\otimes V_{\lambda_3})\\
\otimes Hom(V_{\lambda_0}, V_{\xi}\otimes V_{\lambda_2})\big)
\end{array}
\ar[d]_{i}
}
\\
Hom(V_{\lambda_0}, V_{\lambda_1}\otimes V_{\lambda_2}\otimes V_{\lambda_3})
\ar[r]^{PR}& 
Hom(V_{\lambda_0}, V_{\lambda_1}\otimes V_{\lambda_3}\otimes V_{\lambda_2}), 
}
\end{displaymath}
where $|\lambda_1+\lambda_2|\ge \rho\ge |\lambda_1-\lambda_2|$, $|\lambda_3+\rho|\ge \lambda_0\ge |\lambda_3-\rho|$,
$|\lambda_1+\lambda_3|\ge \xi\ge |\lambda_1-\lambda_3|$, $|\lambda_2+\xi|\ge \lambda_0\ge |\lambda_2-\xi|$, and $i$ is an isomorphism.
\end{prop}

Using the notation, we introduced above, one can write the statements of Proposition \ref{intertwiner1} 
as follows, similarly to (\ref{int1}), (\ref{int2}):
\begin{eqnarray}\label{intfd}
&&(1\otimes PR)\phi_\rho^{\lambda_1\lambda_2}\phi_{\lambda_0}^{\rho\lambda_3}=
\sum_{\xi}B^V_{\rho\xi}\left[\begin{array}{cc}
\lambda_0 & \lambda_1\\
\lambda_2 & \lambda_3
\end{array}
\right]
\phi_\xi^{\lambda_1\lambda_3}\phi_{\lambda_0}^{\xi\lambda_2},\\
&&\label{intfdnew}
\phi_{\rho\lambda_3}^{\lambda_0}\phi_{\lambda_1\lambda_2}^\rho(1\otimes PR)=
\sum_{\xi}B^V_{\xi\rho}\left[\begin{array}{cc}
\lambda_0 & \lambda_1\\
\lambda_2 & \lambda_3
\end{array}
\right]
\phi_{\xi\lambda_2}^{\lambda_0}\phi_{\lambda_1\lambda_3}^\xi.
\end{eqnarray} 
 
It is well known that the relations  (\ref{intfd}), (\ref{intfdnew}) provide the main structure coefficients of the 
braided tensor category $\mathcal{C}_q$ of finite-dimensional representations of $U_q(sl(2))$. Similarly, the relations (\ref{int1}) and (\ref{int2}) allow us to define a generic version of 
$\mathcal{C}_q$ which we denote by $\tilde{\mathcal{C}_q}$. We define objects of $\tilde{\mathcal{C}_q}$ to be infinite sums of modules from the usual category $\mathcal{O}$ of $U_q(sl(2))$, though we still may require that all the weight spaces are finite dimensional. This is sufficient to include the tensor products of Verma modules. Note that $Hom$ spaces between single Verma module and tensor products of Verma modules are finite-dimensional. Then the relations (\ref{int1}), (\ref{int2}) again determine the structure coefficients of the category $\tilde{\mathcal{C}_q}$ at the generic weights. For the integral weights the dimensions of $Hom$ spaces might jump up (in particular, $dim Hom(M_{\nu},M_{\mu}\otimes M_{\lambda})$ might be greater than $1$, see e.g. \cite{chari})
and in order to extend the structure of the category $\tilde{\mathcal{C}_q}$ to the integral weights, one needs a careful study of the analytic continuation of the structure coefficients. In the next subsection we give an explicit realization of the simple objects and $Hom$'s in both categories $\mathcal{C}_q$ and $\tilde{\mathcal{C}_q}$ that will allow us to relate directly their 
structure coefficients.\\


\noindent {\bf {2.3. Polynomial realization for $M_\lambda^c$ and the formula for intertwining operator.}} 
Let's consider two variables $\beta,\zeta$. We claim that the 
space $F_\lambda=\mathbb{C}[\beta]\zeta^\lambda$ carries a structure of $U_q(sl(2))$ module and one can 
identify it with $M_\lambda^c$. Let's introduce $\gamma=\partial^q_\beta$, where $\partial^q_\beta$ is 
a Jackson's q-derivative:
\begin{equation}
\partial^q_\beta f(\beta)=\frac{f(q\beta)-f(q^{-1}\beta)}
{\beta(q-q^{-1})}.
\end{equation}
We denote $v_{m,\lambda}\equiv\beta^m\zeta^\lambda\in\widetilde M_\lambda$. 
These vectors span all $F_\lambda$. Moreover, the following statement holds.\\
\begin{prop} \label{realiz}
Let $\partial_{\beta}, \partial_{\zeta}$ denote usual partial derivatives with respect to $\beta, \gamma$ correspondingly. Then the following identification:
$$
E=q^H\gamma,\quad F=\beta[\zeta\p_{\zeta}-N]q^{-H},\quad H=\zeta\p_{\zeta}-2N,
$$
where $N=\beta\p_{\beta}$ is a number operator (for $a =\zeta\p_{\zeta}-N$, $[a]=\frac{q^a-q^{-a}}{q-q^{-1}}$), 
gives a structure of $U_q(sl(2))$-module on $F_\lambda$, such that $F_\lambda$ is isomorphic to $M_\lambda^c$. 
\end{prop}
\noindent\underline{\bf Proof.\ }
One can find that the action of generators on basis vectors 
$v_{m,\lambda}$ is given by
\begin{eqnarray}
q^{-H}Ev_{m,\lambda}&=&[m]v_{m-1,\lambda},\nonumber\\
Fq^{H}v_{m,\lambda}&=&[\lambda-m]v_{m+1,\lambda},\nonumber\\
q^{\pm H}v_{m,\lambda}&=&q^{\pm(\lambda-2m)}v_{m,\lambda}.
\end{eqnarray}
The resulting module over $\widetilde M_\lambda$ is isomorhpic to $M_\lambda^c$.
Moreover, one can easily get that $v_{m,\lambda}$ corresponds to $[m]!(F^mv_\lambda)^*$, 
where $v_\lambda$ is the highest weight vector in $M_\lambda$, and $\langle v_\lambda^*,v_\lambda\rangle=1$. 
\hfill$\blacksquare$\bigskip

\noindent Next we obtain an explicit form for the intertwining operator 
 $\Phi_{\mu\lambda}^\nu:M_\mu^c\otimes V_\lambda\rightarrow M_\nu^c$ in the realization of Proposition \ref{realiz}. 

Let us denote the coefficients of $\Phi_{\mu\lambda}^\nu$ in such a way:
\begin{equation}
\Phi_{\mu\lambda}^{\nu}(v_{m,\mu}\otimes v_{\ell,\lambda})=
\left(\begin{array}{ccc}
\mu & \lambda & \nu\\
m &\ell & n
\end{array}
\right)
v_{n,\nu},
\end{equation}
where we have $(\lambda-2\ell)+(\mu-2m)=(\nu-2n)$.
We also make the following notation $s=\frac{\lambda+\mu-\nu}{2}$.
At first, we consider the case $\ell=0$.
From the basic property of the intertwiner
\begin{equation}
\t E\Phi_{\mu\lambda}^\nu=q^{2}\Phi_{\mu\lambda}^\nu(q^{-H}\otimes \t E+\t E\otimes 1),\quad {\rm{where}} \quad \t E=q^{-H} E, 
\end{equation}
we get a recurrent relation:
\begin{eqnarray}\label{eact}
\left(\begin{array}{ccc}
\mu & \lambda & \nu\\
m &\ell  &n
\end{array}
\right)
[n]&=&
q^{2m-\mu+2}[\ell]
\left(\begin{array}{ccc}
\mu & \lambda & \nu\\
m&  \ell-1 & n-1
\end{array}
\right)\nonumber\\
&+&
q^2[m]
\left(\begin{array}{ccc}
\mu &\lambda &\nu\\
m-1 &\ell &n-1
\end{array}
\right).
\end{eqnarray}
For $\ell=0$ one obtains:
\begin{equation}
\left(\begin{array}{ccc}
\mu &\lambda &\nu\\
m &0 &n
\end{array}
\right)
[n]=q^2
\left(\begin{array}{ccc}
\mu &\lambda &\nu\\
m-1 & 0  &n-1
\end{array}
\right)
[m].
\end{equation}
We normalize the intertwining operator by the condition:
\begin{equation}
\left(\begin{array}{ccc}
\mu &\lambda &\nu\\
s &0 &0
\end{array}
\right)
=1.
\end{equation}
Therefore,
\begin{equation}
\Phi_{\mu\lambda}^\nu(v_{m,\mu}\otimes v_{0,\lambda})=q^{2(m-n)}
\left[\begin{array}{c}
m \\
n
\end{array}
\right]
v_n^{(\nu)},
\end{equation}
and one can write the explicit formula for $\Phi_{\mu\lambda}^\nu(\cdot\otimes v_{0,\lambda})$
in the polynomial realization. Namely,
\begin{equation}\label{v0}
\Phi_{\mu\lambda}^\nu(v_{m,\mu}\otimes v_{0,\lambda})=\frac{\zeta^{\nu-\mu}q^{2s}}
{[s]!}( \partial^q_\beta )^s v_{m,\mu}.
\end{equation}
In order to obtain the general formula, we use again the basic property of the intertwiner, namely:
\begin{equation}\label{fact}
\widetilde{F}\Phi_{\mu\lambda}^\nu=q^{-2}\Phi_{\mu\lambda}^\nu(1\otimes\widetilde{F}+
\widetilde{F}\otimes q^H), \quad \mbox{{\rm where}} \quad \widetilde{F}=Fq^H.
\end{equation}
Therefore,
\begin{equation}
q^{-2}[\lambda-\ell]\Phi_{\mu\lambda}^\nu(\cdot\otimes v_{\ell+1,\lambda})=-q^{-2}
\Phi_{\mu\lambda}^\nu(\widetilde{F}\cdot\otimes q^{\lambda-2\ell}v_{\ell,\lambda})+
\widetilde{F}\Phi_{\mu\lambda}^\nu(\cdot\otimes v_{\ell,\lambda}).
\end{equation}
Hence,
\begin{eqnarray}
\Phi_{\mu\lambda}^\nu(\cdot\otimes v_{\ell,\lambda})&=&
\frac{q^{2\ell}}{[\lambda-(\ell-1)]\dots[\lambda]}\cdot\nonumber\\
& &\sum_{k=0}^{\ell}\widetilde{F}^{\ell-k}q^{(\lambda -2)k}\frac{\zeta^{\lambda-2s}}{[s]!}
(\partial^q_{\beta})^s(-1)^k \widetilde{F}^k g_k(\ell,q),
\end{eqnarray}
where 
\begin{equation}
g_k(\ell,q)=\sum_{0\leqslant r_1<\dots<r_k\leqslant\ell-1}
q^{-2(r_1+\dots+r_k)}=q^{-k(\ell-1)}
\left[\begin{array}{c}
\ell \\
k
\end{array}
\right]_q.
\end{equation}
Therefore, the following statement holds.

\begin{prop}\label{expint}
Let $\lambda\in \mathbb{Z}_+$. Then 
 the polynomial realization for the operator
\begin{equation}
\Phi_{\mu\lambda}^\nu(v_{\ell,\lambda})\equiv
\Phi_{\mu\lambda}^\nu(\cdot\otimes v_{\ell,\lambda}):
M_\mu^c\rightarrow M_\nu^c,
\end{equation}
where 
$\Phi_{\mu\lambda}^\nu\in Hom(M_\mu^c\otimes V_\lambda,M_\nu^c)$, such that $\frac{\mu+\lambda-\nu}{2}=s$ in the 
case of $\ell=0$, is given by (\ref{v0}) and, if $\ell>0$, the explicit expression is:
\begin{eqnarray}\label{realint}
&&\Phi_{\mu\lambda}^\nu(v_{\ell,\lambda})\cdot=
\frac{q^{2\ell}}{[\lambda]\dots[\lambda-(\ell-1)][s]!}\cdot\nonumber\\
&&\sum_{k=0}^\ell (-1)^kq^{k(\lambda-\ell-1)}
\left[\begin{array}{c}
\ell \\
k
\end{array}
\right]_q
\widetilde{F}^{\ell-k}\zeta^{\lambda-2s}(\partial^q_{\beta})^s\widetilde{F}^k,
\end{eqnarray}
where
\begin{equation}
\t F=\beta\Big[\zeta\frac{\ud}{\ud\zeta}-N\Big].
\end{equation}
In the case of generic $\lambda$, formula (\ref{realint}) gives the polynomial realization for the intertwiner from 
$Hom(M_\mu^c\otimes M^c_\lambda,M_\nu^c)$
\end{prop}

Note, that the operator, we have constructed above, represents an element from $Hom(M_\mu^c\otimes V_\lambda, M_\nu^c)$. 
We note, that in \cite{stv} the intertwiners $Hom(M_{\nu}, M_\mu\otimes V_\lambda)$ were studied in the higher rank case. 

One can show that the intertwiner $Hom(M_\mu^c\otimes V_\lambda,M_\nu^c)$, which we have constructed in Proposition \ref{expint}, 
can be extended in a unique way to another one, from $Hom(M_\mu^c\otimes M^c_\lambda, M_\nu^c)$. 
 
Really, one can allow $\ell$ to take values below $\lambda$ and, therefore, one can 
write down the expression (\ref{eact}) in the case when $\ell=\lambda+1$
 \begin{eqnarray}
\left(\begin{array}{ccc}
\mu & \lambda & \nu\\
m &\lambda+1  &n
\end{array}
\right)
[n]&=&
q^{2m-\mu+2}[\lambda+1]
\left(\begin{array}{ccc}
\mu & \lambda & \nu\\
m&  \lambda & n-1
\end{array}
\right)\nonumber\\
&+&
q^2[m]
\left(\begin{array}{ccc}
\mu &\lambda &\nu\\
m-1 &\lambda+1 &n-1
\end{array}
\right).
\end{eqnarray}
The expression above gives a recurrent relation for the matrix elements 
\begin{eqnarray}\label{coeff}
\left(\begin{array}{ccc}
\mu &\lambda &\nu\\
m &\lambda+1 &n \end{array}
\right)
\end{eqnarray}
and allows to express them by means of the elements $\left(\begin{array}{ccc}
\mu &\lambda &\nu\\
m &\lambda &n 
\end{array}
\right)$
which are the coefficients we already know from the previous calculations for the intertwiner with $M^c_{\lambda}$ reduced to $V_{\lambda}$. 
Once we know the expression for (\ref{coeff}), we can calculate the matrix coefficients of $\Phi_{\mu\lambda}^\nu(v_{\lambda+1,\lambda})$. The coefficients for 
$\Phi_{\mu\lambda}^\nu(v_{\lambda+k,\lambda})$, where $k>1$ can be deduced as before, by means of the action of $\t F$ operator (\ref{fact}). Hence we have the following statement.
\begin{col}\label{col}
There is a unique extension of $\Phi_{\mu\lambda}^{\nu}\in Hom(M^c_{\mu}\otimes V_{\lambda}, M^c_{\nu})$ to the intertwining operator ${\Phi_{\mu\lambda}^{\nu}}'\in Hom(M^c_{\mu}\otimes M^c_{\lambda}, M^c_{\nu})$, such that ${\Phi_{\mu\lambda}^{\nu}}'(\cdot, i_{\lambda} \cdot)=\Phi_{\mu\lambda}^{\nu}(\cdot, \cdot)$, where $i_{\lambda}$ is an inclusion $i_{\lambda}: V_{\lambda}\to M^c_{\lambda}$.
\end{col}

The construction above shows that there exists a continuation of intertwiners from the generic values of weights to the integer values. Therefore, restricting the relation (\ref{int2}) to the subspaces, which in the case of integer $\lambda_i$ $(i=0,1,2,3)$ corresponds to the embedding of the irreducible finite-dimensional modules, we find out that the relation between braiding matrices is the one provided by the Proposition below.

\begin{prop} Let $\lambda_i\in\mathbb{Z}_+$ $(i=0,1,2,3)$.  
There exists a continuation of the elements of the braiding matrix $B_{\rho,\xi}^M$ such that 
\begin{eqnarray}
B^M_{\rho\xi}\left[\begin{array}{cc}
\lambda_0 & \lambda_1\\
\lambda_2 & \lambda_3
\end{array}\right]=
B^V_{\rho\xi}\left[\begin{array}{cc}
\lambda_0 & \lambda_1\\
\lambda_2 & \lambda_3
\end{array}\right],
\end{eqnarray}
where $\rho,\xi\in \mathbb{Z}_+$, such that $\lambda_1+\lambda_2\ge \rho\ge |\lambda_1-\lambda_2|$, $\lambda_3+\rho\ge \lambda_0\ge |\lambda_3-\rho|$,
$\lambda_1+\lambda_3\ge \xi\ge |\lambda_1-\lambda_3|$, $\lambda_2+\xi\ge \lambda_0\ge |\lambda_2-\xi|$.
\end{prop}
\noindent\underline{\bf Proof.\ } In this section we gave explicit construction of the intertwining operator from $Hom(M^c_{\mu}\otimes M^c_{\lambda},M^c_{\nu})$. It is important that the construction holds both in the case of generic points and in the case of integer weights. It makes sense to consider the element $\Phi^{\nu}_{\t\mu\t\lambda}\in Hom(M^c_{\mu+\epsilon_2}\otimes M^c_{\lambda+\epsilon_1},M^c_{\nu+\epsilon_1+\epsilon_2})$, where 
$0<\epsilon_{1,2}<< 1$, $\lambda, \mu \in \mathbb{Z}_+$ and the normalization of this intertwiner is the one from Proposition \ref{expint}.  

Let us consider the identity (\ref{int2})  in the case of the $\epsilon$-regularized intertwiners, such that 
$\t \lambda_i> 0$ (i=0,1,2,3).  We will consider two limits with respect to regularization parameters. At first, we 
evaluate the limit $\t\lambda_{2}\to \lambda_2$, $\t\lambda_{3}\to \lambda_3$. 
Since the intertwining operators exist in this case (see Proposition 2.4. and Corollary 2.1), the limit  of the corresponding braiding matrix elements 
 $B^M_{\rho\xi}\left[\begin{array}{cc}
\t \lambda_0 & \t \lambda_1\\
\lambda_2 & \lambda_3
\end{array}
\right]$
at integer points $\lambda_2$, $\lambda_3$ is finite. Now it possible to restrict the intertwining operators to $V_{\lambda_2}\subset M^c_{\lambda_2}$, $V_{\lambda_3}\subset M^c_{\lambda_3}$. In similar manner one can take the limit $\t \lambda_1\to \lambda_1$ and again, thanks to Proposition 2.4, the limit of the elements of the braiding matrix exists. Then it is possible to restrict interwiners to $V_{\lambda_1}\subset M^c_{\lambda_1}$, i.e. we get braiding relations between operators $\phi^{\nu}_{\mu\lambda}$. Hence, one obtains that the analytical continuation of the elements of the braiding matrix to integer points exists and, when $\lambda_i, \rho,\xi$ satisfy the conditions stated in the proposition, they coincide with the appropriate elements of the braiding matrix $B^V$.  
\hfill$\blacksquare$\bigskip

It is well known that the structure coefficients $B^V_{\rho\xi}\left[\begin{array}{cc}
\lambda_0 & \lambda_1\\
\lambda_2 & \lambda_3
\end{array}\right]$ of the category $\mathcal{C}_q$ can be expressed by quantum $6j$-symbols, which in its turn are given by the balanced basic hypergeometric functions $_4\phi_3$ with the integral values of parameters \cite{KiRe}. Similarly, the structure coefficients $B^M_{\rho\xi}\left[\begin{array}{cc}
\lambda_0 & \lambda_1\\
\lambda_2 & \lambda_3
\end{array}\right]$ of the category $\tilde{\mathcal{C}}_q$ are also given by the functions $_4\phi_3$ with three arbitrary complex parameters corresponding to $\lambda_1, \lambda_2, \lambda_3$ and three other parameters corresponding to $\lambda_0, \rho, \xi$, restricted by the integrality conditions of Proposition 2.1. 
This gives an analytic continuation of $6j$-symbols in 3 out of 6 parameters as well as relations between them.

\section{Braided vertex operator algebra on the lattice of Fock spaces}

\setcounter{equation}{0}
{\bf 3.1. Virasoro algebra: basic facts.}
The Virasoro algebra
\begin{equation}
[L_n,L_m]=(n-m)L_{m+n}+\frac{c}{12}(n^3-n)\delta_{n,-m}
\end{equation}
and its representations has been extensively studied for many years (see e.g. \cite{fb} and references therein). Here we need just basic facts. Let us denote by $M_{h,c}$ and $V_{h,c}$ the Verma module and irreducible module
(with highest weight $h$), correspondingly. 
Throughout the paper we will consider only generic values of $c$.  This means that if we 
parametrize the central charge $c$ in such a way:
 \begin{equation}
c=13-6(\varkappa+\frac{1}{\varkappa}),
\end{equation}
where parameter $\varkappa\in \mathbb{R}\backslash \mathbb{Q}$. Then we have the following proposition (see e.g. \cite{fb}).
\begin{prop}
For generic value of c, Verma module $M_{h,c}$ has a unique singular vector in the case
if $h=h_{m,n}$, where
\begin{equation}\label{virsing}
h_{m,n}=\frac{1}{4}(m^2-1)\varkappa+\frac{1}{4}(n^2-1)\varkappa^{-1}-\frac{1}{2}(mn-1).
\end{equation} 
This singular vector occurs on the level mn, i.e. the value of 
$L_0$ is $h_{m,n}+mn$.
\end{prop}
In the following we will be interested in the modules with
$h=h_{1,n}=\Delta(\lambda)$, where $\lambda=n-1$, $\Delta(\lambda)=-\frac{\lambda}{2}+\frac{\lambda(\lambda+2)}{4\varkappa}$.

\begin{col}Let $c$ be generic and $\lambda \ge 0$, then
$V_{\Delta(\lambda),c}=M_{\Delta(\lambda),c}/M_{\Delta(\lambda)+\lambda+1,c}$, where
$V_{\Delta(\lambda),c}$ is the irreducible Virasoro module with the highest
weight $\Delta(\lambda)$. For $\lambda < 0$ and generic values of $c$ the irreducible module is isomorphic to Verma one, namely, $V_{\Delta(\lambda),c}=M_{\Delta(\lambda),c}$.
\end{col}

\noindent{\bf{3.2. Braided VOA on the lattice of Fock spaces.}} Let us consider the Heisenberg algebra
\begin{equation}
[a_n,a_m]=2\varkappa m\delta_{n+m,0}
\end{equation}
and denote by $F_{\lambda,\varkappa}$ the Fock module associated to this algebra. Namely, 
$F_{\lambda,\varkappa}=S(a_{-1},a_{-2},\dots)\otimes\mathbf{1}_\lambda$, such that
$a_n\mathbf{1}_\lambda=0$ if $n>0$ and $a_0\mathbf{1}_\lambda=\lambda\mathbf{1}_\lambda (\lambda\in\mathbb{C})$, where the elements ${\bf 1}_{\lambda}$ form an additive group (i.e. ${\bf 1}_{\lambda}\cdot{\bf 1}_{\mu}={\bf 1}_{\lambda+\mu}$), isomorphic to $\mathbb{C}$.
It is well known (see e.g. \cite{fb}, \cite{FHL} and references therein) that $F_{0,\varkappa}$ gives rise to the vertex operator algebra, generated by the field
$a(z)=\sum a_nz^{-n-1}$, such that $\deg(a(z))=1$ which has the following operator product expansion (OPE):
\begin{equation}
a(z)a(w)\sim \frac{2\varkappa}{(z-w)^2}.
\end{equation}
We will denote this vertex algebra as $F_{0,\varkappa}(a)$. Moreover, the following is true.

\begin{prop}\label{feiginfuks}
Vertex algebra $F_{0,\varkappa}(a)$ has a vertex operator algebra structure, where the vertex operator, corresponding to the Virasoro 
element, is given by the following formula:
\begin{equation}
L(z)=\frac{1}{4\varkappa}:a(z)^2:+\frac{\varkappa-1}{2\varkappa}a'(z),
\end{equation} 
such that $L(z)=\sum_nL_nz^{-n-2}$ and $L_n$ satisfy Virasoro algebra
relations with the central charge $c=13-6(\varkappa+\frac{1}{\varkappa})$.
\end{prop}

\noindent Now let us consider the following space:
\begin{equation}
\hat{F}_\varkappa=\oplus_{\lambda\in\mathbb{Z}\oplus\mathbb{Z}\varkappa}F_{\lambda,\varkappa}.
\end{equation}
Below we will show that $\hat{F}_\varkappa$ carries a more general structure than VOA, 
namely braided VOA.

Let us consider the following operator:
\begin{equation}\label{xz}
\mathbb{X}(\lambda,z)=\mathbf{1}_\lambda z^{\frac{\lambda a_0}{2\varkappa}}
e^{\big(\frac{\lambda}{2\varkappa}\sum_{n>0}\frac{a_{-n}}{n}z^n\big)}
e^{-\big(\frac{\lambda}{2\varkappa}\sum_{n>0}\frac{a_n}{n}z^{-n}\big)},
\end{equation}
where $\lambda\in\mathbb{Z}\oplus \mathbb{Z}\varkappa$.
It is clear that $\mathbb{X}(\lambda,z)\mathbf{1}_0|_{z=0}=\mathbf{1}_\lambda$.
Moreover, denoting
\begin{equation}
\mathbb{X}_{n_1,\dots,n_k}(\lambda,z)\equiv :a^{(n_1)}(z)\dots a^{(n_k)}(z)\mathbb{X}(\lambda,z):,
\end{equation}
where the symbol ": :'' stands for the Fock space normal ordering and 
$a^{(n)}(z)=\frac{1}{n!}\big(\frac{\ud}{\ud z}\big)^na(z)$, one
can see that
\begin{equation}
\mathbb{X}_{n_1,\dots,n_k}(\lambda,z)\mathbf{1}_0|_{z=0}=a_{-n_1}\dots a_{-n_k}\mathbf{1}_\lambda.
\end{equation}
In such a way, we build the correspondence
\begin{equation}
v\rightarrow Y(v,z)=\sum_{n\in \mathbb{Z}}v_{(n)}z^{-n-1},
\end{equation}
such that $v\in\hat{F}_\varkappa$ and $v_{(n)}\in \mathrm{End}(\hat{F}_\varkappa)$. 

\noindent Let $|z|>|w|$, then
\begin{equation}
\mathbb{X}(\lambda,z)\mathbb{X}(\mu,w)=(z-w)^{\frac{\lambda\mu}{2\varkappa}}
(\mathbb{X}(\lambda+\mu,w)+\dots),
\end{equation}
where dots stand for the terms regular in $(z-w)$.
Let us consider the paths
\begin{eqnarray}
w(t)&=&\frac{1}{2}\big((z+w)+(w-z)e^{\pi it}\big),\nonumber\\
     z(t)&=&\frac{1}{2}\big((z+w)+(z-w)e^{\pi it}\big),
\end{eqnarray}
where $t\in[0,1]$, and let $\mathscr{A}_{z,\omega}$ denote 
the monodromy around these paths, then
\begin{eqnarray}\label{acon}
\mathscr{A}_{z,w}\big(\mathbb{X}(\lambda,z)\mathbb{X}(\mu,w)\big)=
q^{\frac{\lambda\mu}{2}}\mathbb{X}(\mu,w)\mathbb{X}(\lambda,z),
\end{eqnarray}
where $q=e^{\frac{\pi i}{\varkappa}}$. The expression above should be understood in a $weak$ sense, i.e. the analytical continuation 
is performed for the matrix elements of the corresponding operator products. 
Moreover, the matrix elements of operator product expansion (i.e. correlator, see below) $\mathbb{X}(\lambda,z)\mathbb{X}(\mu,w)$ exist in the domain $|z|>|w|$ and the analytical continuation relates it to the matrix elements of operator product expansion $\mathbb{X}(\mu,w)\mathbb{X}(\lambda,z)$, which converge in the domain $|w|>|z|$.

The same property (\ref{acon}) is true if we
take $\mathbb{X}_{n_1,\dots,n_k}(\lambda,z)$ instead of $\mathbb{X}(\lambda,z)$ operators.
The following statement summarizes all the properties of the resulting algebraic object, which is a particular example   
of the $braided$ vertex operator algebra, which we will consider in Section $5$.

\begin{prop}\label{brvoab}
\hfill
\begin{itemize}
\item[1)] There exists a linear correspondence
\begin{equation}
v\rightarrow Y(v,z)=\sum_{n\in\mathbb{Z}}v_{(n)}z^{-n-1}
\end{equation}
such that $v\in\hat{F}_\varkappa$ and $v_{(n)}\in\mathrm{End} \hat{F}_\varkappa$.
\item[2)] Let $|z|>|w|$ and $v_\xi\in F_{\xi,\varkappa}, v_\eta\in F_{\eta,\varkappa}$, where $\xi, \eta \in \mathbb{C}$. Then 
\begin{equation}
\mathscr{A}_{z,w}\big(Y(v_\xi,z)Y(v_\eta,w)\big)=q^{\xi\eta/2}Y(v_\eta,w)Y(v_\xi,z).
\end{equation}
\item[3)] There is a vector $\mathbf{1}=\mathbf{1}_0$, which satisfies
\begin{equation}
Y(\mathbf{1},z)=\mathrm{Id}_{\hat{F}_\varkappa},\qquad Y(v,z)\mathbf{1}|_{z=0}=v
\end{equation}
for any $v\in \hat{F}_\varkappa$.
\item[4)] There exists an element $D\in\mathrm{End}(\hat{F}_\varkappa)$ such that
\begin{equation}
D\mathbf{1}=0,\qquad [D,Y(v,z)]=\frac{\ud}{\ud z}Y(v,z),\qquad \forall v\in \hat{F}_\varkappa.
\end{equation}
\item[5)] There exists an element $\omega\in \hat{F}_\varkappa$ such that
\begin{equation}
Y(\omega,z)=\sum_{n\in\mathbb{Z}}L_nz^{-n-2}
\end{equation}
and $L_n$ satisfy the relations of Virasoro algebra with $L_{-1}=D$.
\end{itemize}
\end{prop}
%
%
We note here, that similar objects were studied in \cite{LD}, where they were called "abelian intertwining algebras".\\ 

\noindent In the following we will consider the subalgebra of the described above braided vertex operator
algebra, related to the subspace $\widetilde{F}_\varkappa=\oplus_{\lambda\in\mathbb{Z}}F_{\lambda,\varkappa}$.
Let's introduce two vertex operators $\mathbb{X}_s^+(z)=\mathbb{X}(-2,z)$ and $\mathbb{X}_s^-(z)=\mathbb{X}(2\varkappa,z)$, which in the physics literature are usually called "screening operators". We have the following statement.

\begin{lemma}
Operators $\mathbb{X}_s^\pm(z)$ have conformal dimension 1, i.e
\begin{equation}
[L_n,\mathbb{X}_s^\pm(w)]=\partial_w\big(w^{n+1}\mathbb{X}_s^\pm(w)\big).
\end{equation}
\end{lemma}

\vspace{5mm}
\noindent\underline{\bf Proof.} It is easy to see that the OPE of $L(z)$ with $\mathbb{X}_s^\pm$ has the following form:
\begin{equation}
L(z)\mathbb{X}_s^\pm(w)\sim \frac{\mathbb{X}_s^\pm(w)}{(z-w)^2}+
\frac{\partial\mathbb{X}_s^\pm(w)}{(z-w)}.
\end{equation}
Then the statement of Lemma follows from the Cauchy theorem.
\hfill$\blacksquare$\bigskip

\noindent The expressions of the form  
\begin{equation}
 \langle v_1^*, Y(u_n,z_n)...Y(u_1,z_1)v_2\rangle,
\end{equation}
where $|z_n|>...>|z_1|>0$,  $v_2, u_1,..., u_n\in \hat {F}_{\varkappa}$ and $v^*_1$ is the element of the dual space 
$\hat {F}^*_{\varkappa}$, are usually called $correlators$.

One can prove (via Wick theorem for the Fock space normal ordering) the following Lemma (see e.g. \cite{cft}), which gives the explicit value for the correlator of operators (\ref{xz}). 
\begin{lemma} Let $z_i\in \mathbb{R}$ $(i=1,...,n)$, such that $z_n>...>z_1>0$ and $\bf{1}_{\nu}^*$ is the highest weight element in the dual Fock space such that $\langle \bf{1}_{\nu}^*, \bf{1}_{\mu}\rangle=\delta_{\nu,\mu}$. Then
\begin{equation}
\langle {\bf{1}}_{\nu}^*,\mathbb{X}(\mu_n,z_n)...\mathbb{X}(\mu_1,z_1) {\bf{1}}_{\mu_0}\rangle=\delta_{\nu, \mu_n+...+\mu_1+\mu_0}\prod_{i>j}(z_i-z_j)^{\frac{\mu_i\mu_j}{2\varkappa}}.
\end{equation}
\end{lemma}
In the following we will be interested in the functions corresponding to the specific type of correlators, namely: 
\begin{eqnarray}
&&\langle {\bf{1}}_{\nu}^*, \mathbb{X}(\lambda_n,z_n)\dots\mathbb{X}^+_s(x_\ell)\dots
\mathbb{X}_s^+(x_1)\dots \mathbb{X}(\lambda_1,z_1) {\bf{1}}_{\lambda_0}
\rangle=\nonumber\\
&&\Psi_{\vec{z}}(x_1,\dots,x_\ell)\delta_{\nu, \lambda_n+...+\lambda_1+\lambda_0-2\ell} , 
\end{eqnarray}
where $\vec{z}$ stands for $(z_1,\dots, z_n)$ and 
 \begin{eqnarray}\label{psi}
\Psi_{\vec{z}}(x_1,\dots,x_\ell)=\prod_{i<j}(x_i-x_j)^{2/\varkappa}
\prod_{j,p}(x_i-z_p)^{-\lambda_p/\varkappa}
\prod_{p<q}(z_p-z_q)^{\frac{\lambda_p\lambda_q}{2\varkappa}}.
\end{eqnarray}
This will motivate our constructions in the next section, where we will consider the local system generated by the branches of this function. 

\noindent In the next subsection, we will study the relation between the Fock spaces and Virasoro
modules.\\

\noindent{\bf{3.3. Irreducible highest weight Virasoro modules and the screening charge.}} In previous subsection, we introduced Fock modules 
$F_{\lambda,\varkappa}$ and also screening operators $\mathbb{X}_s^\pm(z)$. 
Next Proposition will explain how to construct the so-called "screening charge'', associated with 
$\mathbb{X}_s^-(z)$, which acts on $F_{\lambda,\varkappa}$.

\begin{prop}\label{qminus}
Let $Q^-$ be the operator $Q^-:F_{\lambda,\varkappa}\rightarrow F_{\lambda+2\varkappa,\varkappa}$,
where $\lambda\in\mathbb{Z}$, which is defined by means of the following
formula:
\begin{equation}
Q^-v=\oint_{C_{z_2}}\frac{\ud z}{2\pi i}\mathbb{X}_s^-(z)v,
\quad v\in F_\lambda , 
\end{equation}
then
\begin{itemize}
\item[(i)] action of  $Q^-$ is well defined on $F_{\lambda,\varkappa}$, i.e. operator 
$\mathbb{X}_s^-(z_1)$ has only integer powers in the OPE with $Y(v,z_2)$;
\item[(ii)] the operator $Q^-$ commutes with the action of Virasoro operators: \newline
$[Q^-,L_n]=0$.
\end{itemize}
\end{prop}

\vspace{5mm}
\noindent\underline{\bf Proof.}
(i) follows from the fact that 
\begin{equation}
\mathbb{X}_s^-(z)\mathbb{X}(\lambda,w)=(z-w)^\lambda\big(\mathbb{X}(\lambda+2\varkappa,w)+\dots\big),
\end{equation}
where dots stand for the terms with regular powers in $(z-w)$.

(ii) comes from the following calculation:
\begin{equation}
[L_n,Q^-]=\int_C\frac{\ud w}{2\pi i}\int_{C_{w}}\frac{\ud z z^{n+1}}{2\pi i}
L(z)\mathbb{X}_s^-(w)=
\int_C\ud w\partial_w\big(w^{n+1}\mathbb{X}_s^-(w)\big).
\end{equation}
\hfill$\blacksquare$\bigskip

The kernel of $Q^-$ gives an irreducible representation of Virasoro algebra, namely, we have the 
following statement.
\begin{prop}
The kernel of the operator $Q^-$ acting on $F_{\lambda,\varkappa}$, where $\lambda \in \mathbb{Z}$, is given by the following expressions:  
\begin{eqnarray}
\ker Q^-|_{F_{\lambda,\varkappa}}&=&V_{\Delta(\lambda),\varkappa},\qquad \lambda\geqslant 0,\quad 
\Delta(\lambda)=-\frac{\lambda}{2}+\frac{\lambda(\lambda+2)}{4\varkappa},\nonumber\\
\ker Q^-|_{F_{\lambda,\varkappa}}&=&0,\qquad\qquad\lambda<0.
\end{eqnarray}
\end{prop}

\noindent\underline{\bf Proof.}
From the previous Proposition we already know that $Q^-$ maps $F_\lambda$ to $F_{\lambda '}$,
where $\lambda'=\lambda+2\varkappa$. Therefore,
$\Delta(\lambda')=\frac{\lambda+2}{2}+\frac{\lambda(\lambda+2)}{4\varkappa}$.
Let's write the Virasoro characters for $F_\lambda$ and $F_\lambda '$:
\begin{eqnarray}
\mathrm{ch } F_{\lambda,\varkappa}&=&q^{\Delta(\lambda)}\prod_r(1-q^r)^{-1},\nonumber\\
\mathrm{ch } F_{\lambda ',\varkappa}&=&q^{\Delta(\lambda)+\lambda+1}\prod_r(1-q^r)^{-1}.
\end{eqnarray}
Let $\lambda\geqslant 0$, then
\begin{equation}
\mathrm{ch} F_{\lambda,\varkappa}-\mathrm{ch} F_{\lambda',\varkappa}=q^{\Delta(\lambda)}(1-q^{\lambda+1})\prod_r(1-q^r)^{-1}=\mathrm{ch}V_{\Delta(\lambda),\varkappa}.
\end{equation}
So, in order to prove the Proposition we need to show that $ \mathrm{Im}\ Q^-=F_{\lambda,\varkappa}$.
However it is easy to see that the vector, corresponding to the operator 
$:a^{\lambda+1}(z)\mathbb{X}(\lambda,z):$, is mapped by $Q^-$ to the highest
weight state in $F_{\lambda'}$. Therefore, we have proved the statement for $\lambda\geqslant 0$.

In the case of $\lambda<0$ it is not hard to see that the highest weight from $F_{\lambda,\varkappa}$ is not 
mapped to 0 by $Q^-$. Since  $Q^-$ commutes with $L_n$ we find that $\ker Q^-|_{F_{\lambda,\varkappa}}=0$, where $\lambda<0$.
\hfill$\blacksquare$\bigskip
\begin{col}
The space $F_{\lambda,\varkappa}$ gives 
a realization for the dual Verma module of the Virasoro algebra, namely, we have the following exact sequence for $\lambda \ge0$ (cf. (\ref{ex}) ):
\begin{equation}
0\to V_{\Delta(\lambda),\varkappa}\to F_{\lambda,\varkappa}\to V_{\Delta(-\lambda-2),\varkappa}\to 0
\end{equation}
and $F_{\lambda,\varkappa}\cong V_{\Delta(\lambda),\varkappa}$ for $\lambda\in \mathbb{Z}_{\le -1}$.
\end{col}

\section{$\mathbf{U_q(sl(2))}$ and the homology of local systems}
\setcounter{equation}{0}
\noindent{\bf{4.1. Local systems on configuration space and homology.}} In this section, using the combination of approaches of \cite{gs}, \cite{fw}, \cite{varchenko},  we consider geometric versions of the objects we introduced in Section 2, namely finite dimensional irreducible modules, Verma modules and intertwiners, or, in other words, some part of the braided tensor category associated to $U_q(sl(2))$. 
Here we give the sketch of results, considering only the results which we will need for our further constructions. For more details one should consult \cite{fw} and \cite{varchenko}. 
 
The key object for us will be the function $\Psi_{\vec{z}}(x_1,...,x_{\ell})$ (\ref{psi}), which gives the value to the correlator with screening operators $\mathbb{X}^{+}_s$.

We consider the following data: $\lambda_1,\dots,\lambda_n\in\mathbb{C}$, $z_1,\dots,z_n\in\mathbb{C},\quad z_i\ne z_j,\quad z_i\ne 0$, $\varkappa\in\mathbb{R}$$ \backslash $$\mathbb{Q}$. Let $\lambda=\lambda_1+\dots+\lambda_n$, $q=e^{\frac{\pi i}{\varkappa}}$. 
Let us denote by $\mathcal{H}$ the set of hyperplanes 
\begin{eqnarray*}
&&x_i=x_j, \quad i,j=1,\dots,\ell,\\
&&x_i=z_k, \quad i=1,\dots,\ell, \quad k=1,\dots,n.
\end{eqnarray*}
in $\mathbb{C}^{\ell}=\{(x_1,\dots,x_\ell)\}$. Now we consider the 1-dimensional 
local system $\mathscr{S}$ (see e.g. \cite{efk}) over $\mathbb{C}^\ell\backslash\mathcal{H}$
such that its flat sections are $s(x)=\alpha$(univalent branch of $\Psi_{\vec{z}}(x)$) such that
$\alpha\in\mathbb{C}$ and (see \ref{psi})
\begin{eqnarray}
\Psi_{\vec{z}}(x_1,\dots,x_\ell)=\prod_{i<j}(x_i-x_j)^{2/\varkappa}
\prod_{j,p}(x_i-z_p)^{-\lambda_p/\varkappa}
\prod_{p<q}(z_p-z_q)^{\frac{\lambda_p\lambda_q}{2\varkappa}},
\end{eqnarray}
where $\vec{z}$ stands for $(z_1,\dots, z_n)$. 

Now we give a prescription how one should choose a section of the corresponding local system.
For any pair of indices $i,j$ we define
\begin{equation}
Br\big((x_i-x_j)^\rho\big)=
\left\{\begin{array}{cc}
e^{\rho\log(x_i-x_j)}& \mathrm{Re}\  x_i>\mathrm{Re}\  x_j\\
e^{\rho\log(x_j-x_i)}&  \mathrm{Re}\ x_i<\mathrm{Re}\  x_j
\end{array}
\right.
,
\end{equation}
where $\log(x)$ is the principal branch of the logarithm defined for $\mathrm{Re}\ x>0$
by the condition $\log x\in \mathbb{R}_{+}$ for $x>0$. Similarly, we define
\begin{equation}\label{branch}
\mathrm{Br}(\Psi)=\prod\mathrm{Br}(x_i-x_j)^{2/\varkappa}
\prod\mathrm{Br}(x_i-z_k)^{-\lambda_k/\varkappa}
\prod\mathrm{Br}(z_i-z_j)^{\frac{\lambda_i\lambda_j}{2\varkappa}}.
\end{equation}
One can see that if $z_i-z_j\notin i\mathbb{R}$ and a region $D\subset\mathbb{C}^\ell$
 is such that $x_i-x_j\notin i\mathbb{R}$ on $D$, then $\mathrm{Br}(\Psi)$ is
a section of $\mathscr{S}$ over $D$.

One can define twisted $\ell$-cells associated with this local system: it is a pair $(\Delta^{\ell},s)$, where 
$\Delta^{\ell}\subset \mathbb{C}^{\ell}\backslash \mathcal{H}$ is a singular $\ell$-cell and $s$ is the univalent branch of $\Psi$. In such a way one can define a boundary operator on twisted $\ell$-cell 
$\partial_{\ell}$ by the formula $\partial_{\ell}(\Delta^{\ell},s)=(\partial\Delta^{\ell},s |_{\partial\Delta^{\ell}})$. In this  section our main objects will be homology groups $H_\ell(\mathbb{C}^\ell\backslash\mathcal{H},\mathscr{S})$=$\ker\partial_{\ell}/Im\partial_{\ell+1}$ modulo the permutation of the coordinates.
Namely, there is a natural action of the permutation group $\Sigma_\ell$ on a set of coordinates 
$x_1\dots x_\ell$ such that this action preserves $\mathcal{H}$ and $\mathscr{S}$.

So, let us denote the antisymmetric parts of the corresponding homology groups as follows:
\begin{equation}\label{homology}
H^{-\Sigma}_\ell(z_1,\dots,z_n,\lambda_1,\dots,\lambda_n)=
H_\ell^{-\Sigma}(\mathbb{C}^\ell\backslash\mathcal{H},\mathscr{S}).
\end{equation}

\noindent{\bf{4.2. Gomez-Sierra contours, Verma modules and irreducible modules.}} 
In this subsection we introduce an additional construction, which gives the geometric description of the tensor product of Verma modules. We consider the twisted chains of specific kind, namely loops, over which $x_i$ are running, emanating from some fixed point $a\in \mathbb{C}$ around points $z_1,\dots, z_n$ (see e.g. \cite{fw}, \cite{varchenko}, \cite{varchenko2}). For calculations, it is useful to move the reference point $a$ to infinity. 
In such a way, these loops will be transformed to the infinite chains with boundaries involving the point at infinity, see e.g. the chains $G^\ell_{r_1r_2\dots, r_n}(z_1,\dots,z_n,\lambda_1,\dots,\lambda_n)$ on Fig. 1 (see \cite{gs}) together with the associated branch of $\Psi_{\vec{z}}$ chosen in accordance with (\ref{branch}). 
At first they were introduced by Gomez and Sierra in \cite{gs}, therefore in the following we will refer to these loops as Gomez-Sierra contours. These chains also appeared in \cite{fkv} as dual to certain relative homology 
cycles. Note that the authors of \cite{gs} view these chains as integration contours (without caring about convergence) for some  function $f=\Psi_{\vec{z}}(x_1,...,x_{\ell})A(\vec{z},\vec{x})$, where $A(\vec{z},\vec{x})$ is a meromorphic function of $(z_k-z_l)$, $(x_i-x_j)$, $(z_k-x_j)$. 

               
\begin{figure}[hbt] 
\begin{minipage}[c]{0.50\textwidth}
\centering           
\includegraphics[width=0.7\textwidth]{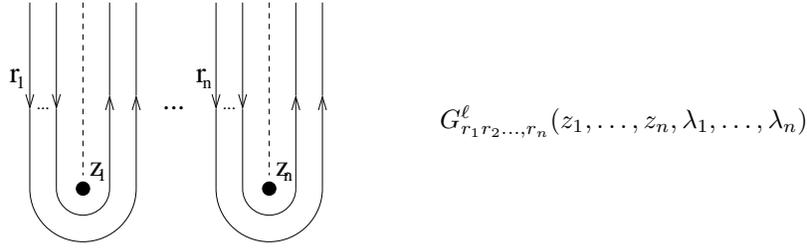}
\end{minipage}%
\begin{minipage}[c]{0.5\textwidth}
\centering
\begin{tabular}{c}
$G^\ell_{r_1r_2\dots, r_n}(z_1,\dots,z_n,
\lambda_1,\dots,\lambda_n)$
\end{tabular}
\end{minipage}
\caption{Basic chains}                
\end{figure}     

Each of the contours on the Fig. 1 goes counter-clockwise 
along the cut (provided by the branches of $\Psi$) and around one of the points $z_i$ in such away that $r_1+...+r_n=\ell$. We can treat them in a formal way as relative cycles with respect to the set of hyperplanes corresponding to the point at infinity ($x_i=\infty$). 
Let us consider the space $\mathcal{S}_\ell(z_1,\dots,z_n;\lambda_1,\dots,\lambda_n)$, spanned by these cycles modulo the homotopy  
trivial ones. We will be interested in the antisymmetric part of it with respect to the action of permutation group $\Sigma$. In the following we denote such space as 
\begin{eqnarray}
\mathcal{S}_\ell^{-\Sigma}(z_1,\dots,z_n;\lambda_1,\dots,\lambda_n).
\end{eqnarray}
Let us denote by $M_{\lambda_1}\otimes\dots\otimes M_{\lambda_n}[\lambda-2\ell]$ 
the subspace of weight $\lambda-2\ell$ (where $\lambda=\sum^n_{i=1}\lambda_i$) of the tensor product of the corresponding Verma modules of $U_q(sl(2))$.
There is a one-to-one a map: 
\begin{eqnarray}
\xymatrix{
M_{\lambda_1}\otimes\dots\otimes M_{\lambda_n}[\lambda-2\ell]
\ar[r]^{\varphi_{\vec{z}}} &
\mathcal{S}^{-\Sigma}_\ell(z_1,\dots,z_n;\lambda_1,\dots,\lambda_\ell), \\
}
\end{eqnarray}  
such that for $k_1+\dots+k_n=\ell$
\begin{equation}
\varphi_{\vec{z}}:\ F^{k_1}v_{\lambda_1}\otimes\dots\otimes F^{k_n}v_{\lambda_n}
\longmapsto
G^{\ell}_{k_1,\dots,k_n}(z_1,\dots,z_n;\lambda_1,\dots,\lambda_\ell).
\end{equation}
Therefore, there is one-to-one map between 
$\oplus_\ell \mathcal{S}_\ell^{-\Sigma}(z_1,\dots,z_n;\lambda_1,\dots,\lambda_n)$ and 
$M_{\lambda_1}\otimes\dots\otimes M_{\lambda_n}$.

Moreover, the geometric realization of the action of $F$-generator is given by the statement below. 
\begin{prop}\label{coprod}
The following diagram is commutative:
\begin{eqnarray}
\xymatrix{
{\begin{array}{c}
M_{\lambda_1}\otimes\dots\otimes M_{\lambda_i}\otimes\\
M_{\lambda_{i+1}}\otimes...\otimes M_{\lambda_n}[\lambda-2\ell]
\end{array}}
\ar[d]^{\Delta_{i,i+1}(F)} \ar[r]^{\varphi_{\vec{z}}} &
\mathcal{S}^{-\Sigma}_\ell(z_1,\dots,z_n;\lambda_1,\dots,\lambda_n) \ar[d]^{\hat{\Delta}_{i,i+1}(F)}\\
{\begin{array}{c}
M_{\lambda_1}\otimes\dots\otimes M_{\lambda_i}\otimes\\
M_{\lambda_{i+1}}\otimes ...\otimes M_{\lambda_n}[\lambda-2\ell-2]
\end{array}}
\ar[r]^{\varphi_{\vec{z}}} &
\mathcal{S}^{-\Sigma}_{\ell+1}(z_1,\dots,z_n;\lambda_1,\dots,\lambda_n).
}
\end{eqnarray}
where the map $\hat{\Delta}_{i,i+1}(F)$ is represented by the Fig. 2.

\begin{figure}[hbt]                                                       
\centering           
\includegraphics[width=1.\textwidth]{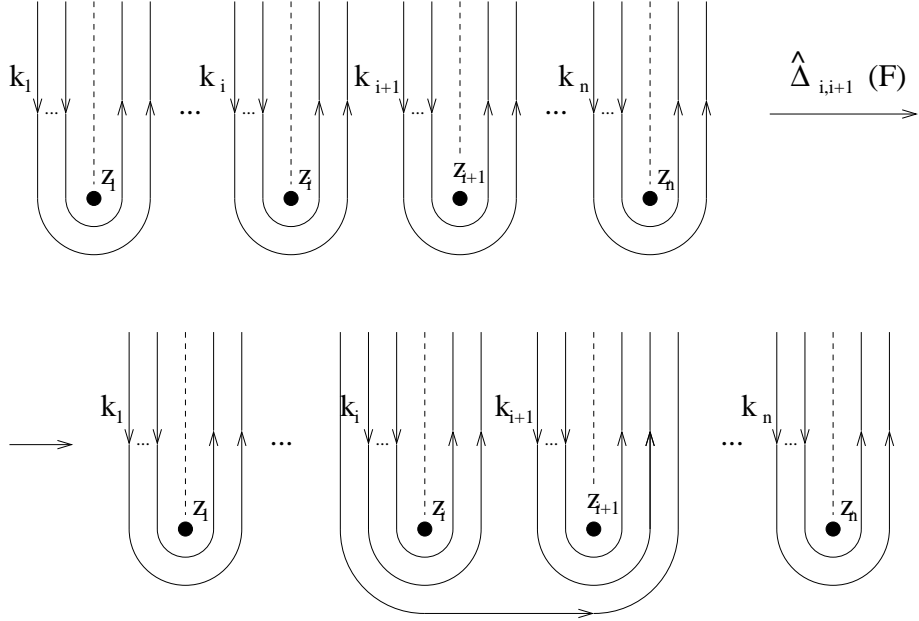}
\caption{The action of coproduct}                
\end{figure}     

\end{prop}

\noindent\underline{\bf Proof.}
In order to prove the Proposition it is enough to do it in the case of two points $z_1,z_2$. 
Namely,

\begin{figure}[hbt] 
\centering           
\includegraphics[width=0.8\textwidth]{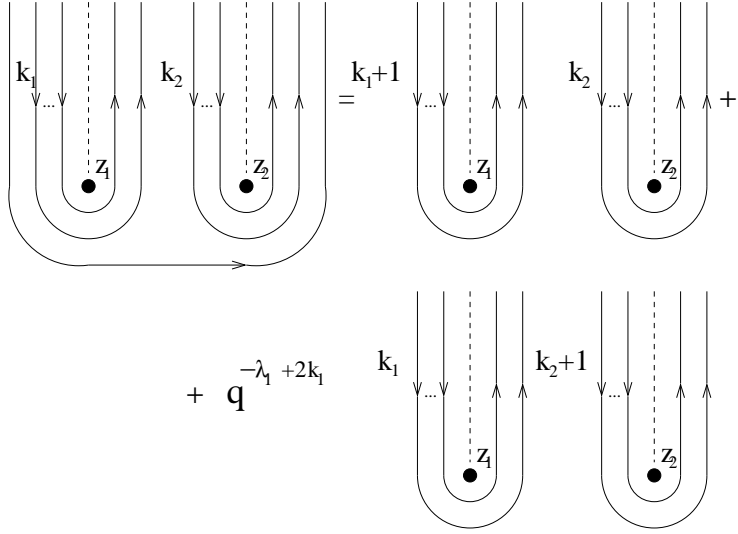}
\caption{The action of coproduct for two points}                
\end{figure}     
\begin{eqnarray}
&&\hat{\Delta}_{1,2}(F)G_{k_1,k_2}^\ell(z_1,z_2;\lambda_1,\lambda_2)=\nonumber\\
&&G_{k_1+1,k_2}^{\ell+1}(z_1,z_2;\lambda_1,\lambda_2)+
q^{-\lambda_1-2k_1}G_{k_1,k_2+1}^{\ell+1}(z_1,z_2;\lambda_1,\lambda_2)=\nonumber\\
&&\varphi_{\vec{z}}\big((F\otimes 1+q^{-H}\otimes F)F^{k_1}v_{\lambda_1}\otimes F^{k_2}v_{\lambda_2}\big)=\nonumber\\
&&\varphi_{\vec{z}}\big(\Delta F(F^{k_1}v_{\lambda_1}\otimes F^{k_2}v_{\lambda_2}\big),
\end{eqnarray}
where the $q$-factor in the second line comes from the structure of branches of $\Psi$, see also Fig. 3.  Therefore, Proposition is proven. 
\hfill$\blacksquare$\bigskip

We can also give the geometric meaning to the action of $F$ on the whole tensor product of Verma modules. One just needs to consider the contour, which embraces not only two (as it was in the Proposition \ref{coprod}), but all $n$ families of Gomez-Sierra contours. 

\noindent The last statement in this subsection gives the geometric meaning to the $R$-matrix.

\begin{prop}\label{rmatrix}

Let $\mathrm{Re}\ z_1<\dots< \mathrm{Re}\ z_i<\mathrm{Re}\ z_{i+1}\dots<\mathrm{Re}\ z_n$. Then the following diagram
 is commutative:
\begin{eqnarray}
\xymatrix{
{M_{\lambda_1}\otimes\dots\otimes M_{\lambda_n}[\lambda-2\ell]}
        \ar[d]^{\check{R}_{i,i+1}} \ar[r]^{\qquad\qquad\varphi_{\vec{z}}}&
\mathcal{S}^{-\Sigma}_\ell(\vec{z},\vec{\lambda})
       \ar[d]^{\mathscr{A}_{i,i+1}}\\
{M_{\lambda_1}\otimes\dots\otimes M_{\lambda_n}[\lambda-2\ell]}
          \ar[r]^{\qquad\qquad\varphi_{s_i(\vec{z})}\qquad} &
\mathcal{S}^{-\Sigma}_\ell\big(s_i(\vec{z}),s_i(\vec{\lambda})\big),
}
\end{eqnarray}
where $s_i(\vec{z})=(z_1,\dots,z_{i+1},z_i,\dots,\dots,z_n),
s_i(\vec{\lambda})=(\lambda_1,\dots,\lambda_{i+1},\lambda_i,\dots,\lambda_n)$.
Here $\mathscr{A}_{i,i+1}$ is the monodromy operator along the path:
\begin{eqnarray}
z_i(t)&=&\frac{1}{2}\big((z_i+z_{i+1})+((z_i-z_{i+1})e^{\pi it}\big),\nonumber\\
z_{i+1}(t)&=&\frac{1}{2}\big((z_i+z_{i+1})+((z_{i+1}-z_i)e^{\pi it}\big),
\end{eqnarray}

\begin{eqnarray}
z_1\ \bullet\  \dots \quad z_i\ 
\xymatrix{
\bullet \ar@/_1pc/[r] &
\bullet   \ar@/_1pc/[l]}
z_{i+1} \quad \dots \ \bullet\ z_n.
\end{eqnarray}

\end{prop}

\vspace{5mm}
%
%
\noindent\underline{\bf Proof.} The proof follows the same steps as in \cite{rrr}. 
In order to prove this we note that is enough to prove it just for two points. In other words, we consider $\mathscr{A}_{12}G_{n_1,n_2}(z_1,z_2,\lambda_1,\lambda_2)$ and reexpress it as follows:
\begin{eqnarray}
&&\mathscr{A}_{12}G_{n_1,n_2}(z_1,z_2,\lambda_1,\lambda_2)=\nonumber\\
&&q^{\frac{\lambda_1(\lambda_2-2n_2)}{2}}\sum_{k=0}^{n_1}\mathrm{C}_{\lambda_1}^{n_1}(k)\big(\hat{\Delta}(F)\big)^{n_1-k}
G_{n_2+k,0}(z_2,z_1,\lambda_2,\lambda_1).
\end{eqnarray}

\begin{figure}[hbt] 
\centering           
\includegraphics[width=0.9\textwidth]{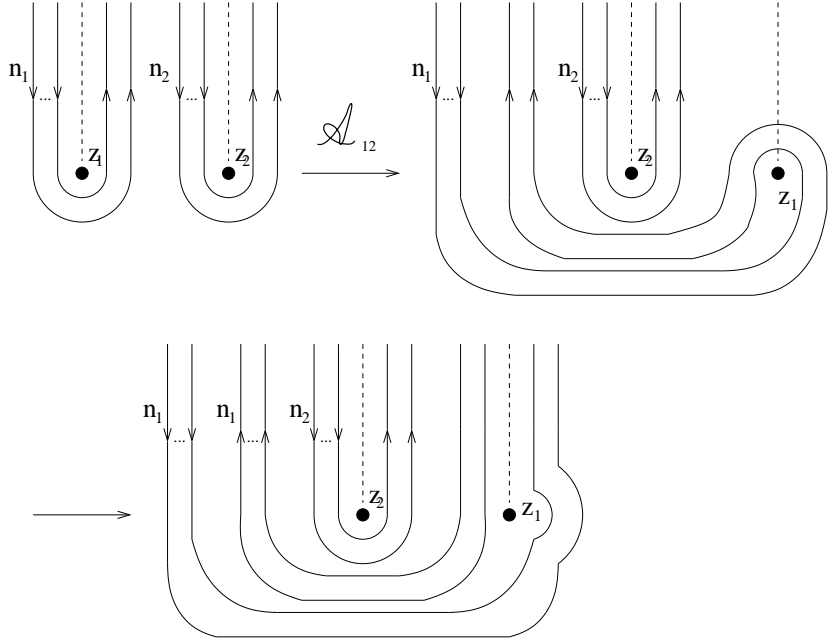}
\caption{Action of the monodromy operator}                
\end{figure}     

In order to determine the coefficients $\mathrm{C}_{\lambda_1}^{n_1}(k)$, we use induction. Suppose, we had $n+1$ contours over $z_1$.

Then the expression in this case can be decomposed as follows:
\begin{eqnarray}
&&\mathscr{A}_{12}G_{n_1+1,n_2}(z_1,z_2,\lambda_1,\lambda_2)=\nonumber\\
&&q^{\frac{\lambda_1(\lambda_2-2n_2)}{2}}\sum_{k=0}^{n_1}\mathrm{C}_{\lambda_1}^{n_1}(k)\big(\hat{\Delta}(F)\big)^{n_1-k+1}
G_{n_2+k,0}(z_2,z_1,\lambda_2,\lambda_1)-\nonumber\\
&&q^{-2\lambda_1+2n_1)}\sum_{k=0}^{n_1}\mathrm{C}_{\lambda_1}^{n_1}(k)\big(\hat{\Delta}(F)\big)^{n_1-k}
G_{n_2+k+1,0}(z_2,z_1,\lambda_2,\lambda_1).
\end{eqnarray}
Therefore, the recurrent relation is:
\begin{equation}
\mathrm{C}_{\lambda_1}^{n_1+1}(k)=\mathrm{C}_{\lambda_1}^{n_1}(k)-q^{-2\lambda_1+2n_1}\mathrm{C}_{\lambda_1}^{n_1}(k-1).
\end{equation}
Let's define the q-deformation of binomial coefficient in the following way:
\begin{equation}
\binom{n}{k}_q\equiv\frac{(n)_q!}{(k)_q!(n-k)_q!},
\end{equation}
where $(n)_q=\frac{q^n-1}{q-1}$ and $(n)_q!=(n)_q(n-1)_q\dots 1$. 
For the binomial coefficients we have the following relation
\begin{equation}
\binom{n+1}{k}_{q^2}=\binom{n}{k}_{q^2}+q^{2n-2k+2}\binom{n}{k-1}_{q^2}.
\end{equation}
Therefore, the expression for $\mathrm{C}_{\lambda}^n(k)$ is the following one:
\begin{equation}
\mathrm{C}_{\lambda}^n(k)=(-1)^kq^{-2k\lambda+k(k-1)}\binom{n}{k}_{q^2}.
\end{equation}
On the other hand,
\begin{eqnarray}
&&\big(\Delta(F)\big)^{n-k}=\sum_{r=0}^{n-k}\binom{n-k}{r}_{q^2}(F^{n-k-r}\otimes F^r)(q^{-rH}\otimes1)=\nonumber\\
&&=\sum_{r=0}^{n-k}q^{-2r(n-k-r)}\binom{n-k}{r}_{q^2}(q^{-rH}\otimes1)(F^{n-k-r}\otimes F^r).
\end{eqnarray}
Therefore,
\begin{eqnarray}
&&\mathscr{A}_{12}G_{n_1,n_2}(z_1,z_2,\lambda_1,\lambda_2)=\nonumber\\
&&q^{\frac{\lambda_1(\lambda_2-2n_2)}{2}}
\sum_{k=0}^{n_1}\sum_{r=0}^{n_1-k}(-1)^kq^{2k\lambda_1+k(k-1)}\binom{n_1}{k}_{q^2}\binom{n_1-k}{r}_{q^2}\nonumber\\
&&\times q^{-2r(n-k-r)}(q^{-rH}\otimes 1)G_{n_1+n_2-r,r}(z_2,z_1;\lambda_2,\lambda_1).
\end{eqnarray}
Using the fact that 
\begin{equation}
\binom{n_1}{k}_{q^2}\binom{n_1-k}{r}_{q^2}=\binom{n_1}{r}_{q^2}\binom{n_1-r}{k}_{q^2}
\end{equation}
and changing the summation, we rewrite the expression above as
\begin{eqnarray}
&&q^{\frac{\lambda_1(\lambda_2-2n_2)}{2}}\sum_{k=0}^{n_1}\sum_{r=0}^{n_1-r}
(-1)^kq^{2k\lambda_1+k(k-1)}\binom{n_1}{r}_{q^2}\binom{n_1-r}{k}_{q^2}\nonumber\\
&&\times q^{-2r(n-k-r)}(q^{-rH}\otimes 1)G_{n_2+n_1-r,r}(\lambda_2,\lambda_1;z_2,z_1).
\end{eqnarray}
Next, we use the formula
\begin{equation}
(1-z)(1-q^2z)\dots(1-q^{2n-2}z)=
\sum_{k=0}^{n}(-1)^k\binom{n}{k}_{q^2}q^{k(k-1)}z^k.
\end{equation}
Let us take $n=n_1-r$ and $z=q^{-2\lambda_1+2r}$. Therefore,
\begin{eqnarray}
&&\big(1-q^{-2\lambda_1+2r}\big)\dots\big(1-q^{2(n_1-r)-2}q^{-2\lambda_1+2r}\big)=\nonumber\\
&&\sum_{k=0}^n(-1)^k\binom{n_1-r}{k}_{q^2}q^{k(k-1)}q^{-2k\lambda_1+2rk}.
\end{eqnarray}
We rewrite the expression above as follows:
\begin{eqnarray}
q^{\frac{\lambda_1(\lambda_2-2n_2)}{2}}\sum_{r=0}^{n_1}\prod_{\ell=1}^{n_1-r}
\big(1-q^{-2\lambda+2n_1-2\ell}\big)q^{-2r(n-r)}\binom{n_1}{r}_{q^2}\times\nonumber\\
(q^{-rH}\otimes 1)G_{n_2+n_1-r,r}(\lambda_2,\lambda_1;z_2,z_1).
\end{eqnarray}
Therefore, the final result is:
\begin{eqnarray}
&&q^{\frac{\lambda_1(\lambda_2-2n_2)}{2}}\sum_{\ell=0}^{n_1}\prod_{s=1}^{\ell}
\big(1-q^{-2\lambda_1+2n_1-2s}\big)q^{-2(\ell-n_1)\ell}\binom{n_1}{\ell}_{q^2}\times\nonumber\\
&&(q^{-(n_1-\ell)H}\otimes 1)G_{n_2+\ell,n_1-\ell}(\lambda_2,\lambda_1;z_2,z_1)=\nonumber\\
&&q^{\frac{\lambda_1(\lambda_2-2n_2)}{2}}\sum_{\ell=0}^{n_1}\prod_{s=1}^{\ell}
\big(1-q^{-2\lambda+2n_1-2s}\big)q^{(-n_1-\ell)(\lambda_2-2n_2)}\times \nonumber\\
&& \binom{n_1}{\ell}_{q^2}G_{n_2+\ell,n_1-\ell}(z_2,z_1;\lambda_2,\lambda_1).
\end{eqnarray}
Now we need another formula:
\begin{equation}
E^\ell F^{n_1}v_{\lambda_1}=\prod_{s=1}^\ell[n_1-s+1][\lambda_1-n_1+s]F^{n-\ell}v_{\lambda_1}.
\end{equation}
We see that 
\begin{eqnarray}
&&\prod_{s=1}^{\ell}\big(1-q^{-2\lambda_1+2n_1-2s}\big)=(q-q^{-1})^\ell\prod_{s=1}^{\ell}q^{-\lambda_1+n_1-s}
\prod_{s=1}^\ell[\lambda_1-n_1+s],\nonumber\\
&&\prod_{s=1}^\ell q^{-\lambda_1+n_1-s}=q^{\ell(-\lambda_1+n_1)}q^{\frac{-\ell(\ell-1)}{2}},\nonumber\\
&&\binom{n_1}{\ell}_{q^2}=\frac{(n_1)_{q^2}!}{(\ell)_{q^2}!(n_1-\ell)_{q^2}!}=
\frac{q^{\frac{n_1(n_1+1)}{2}}}
{q^{\frac{\ell(\ell+1)}{2}}q^{\frac{(n_1-\ell)(n_1-\ell+1)}{2}}}\frac{[n]!}{[\ell]![n_1-\ell]!}=\nonumber\\
&&\frac{1}{[\ell]!}\prod_{s=1}^\ell[n_1-s+1]q^{-\ell^2+n_1\ell}.
\end{eqnarray}
Now let's collect all $q$-factors:
\begin{eqnarray}
&&q^{\frac{\lambda_1(\lambda_2-2n_2)}{2}}q^{\ell(-\lambda_1+n_1)}q^{\frac{-\ell(\ell+1)}{2}}
q^{-\ell^2}q^{-(n_1-\ell)(\lambda_2-2n_2)}=\nonumber\\
&&q^{\frac{(\lambda_1-2n_1+2\ell)(\lambda_2-2n_2-2\ell)}{2}}q^{\frac{\ell(\ell-1)}{2}}.
\end{eqnarray}
Hence, we have
\begin{eqnarray}
&&\mathscr{A}_{12}\big(G_{n_1,n_2}(z_1,z_2;\lambda_1,\lambda_2)\big)=\nonumber\\
&&=\varphi_{z_1,z_2}\bigg(\Big(\sum_{\ell\geqslant 0}q^{\frac{\ell(\ell-1)}{2}}
\frac{(q-q^{-1})^\ell}{[\ell]!}F^\ell\otimes E^\ell\Big)F^{n_2}v_{\lambda_2}\otimes F^{n_1}v_{\lambda_1}\bigg)=\nonumber\\
&&=\varphi_{z_1,z_2}\Big(\check{R}\cdot(F^{n_2}v_{\lambda_2}\otimes F^{n_1}v_{\lambda_1})\Big).
\end{eqnarray}
This finishes the proof. 
\hfill$\blacksquare$\bigskip

\noindent{\bf{4.3. Singular vectors and intertwiners.}} In Proposition \ref{coprod} we defined the action of the $F$-generator and the coproduct in a geometric 
way. However, the action of the $E$-generator is implicitly defined in the following statement.

\begin{prop}\label{cycsing}
Let $\mathrm{Re}\ z_1<\dots< \mathrm{Re}\ z_n$. There is a natural isomorphism between the spaces of singular vectors from
$M_{\lambda_1}\otimes\dots\otimes M_{\lambda_n}[\lambda-2\ell]$ and cycles
$Z_\ell^{-\Sigma}(z_1,\dots,z_n;\lambda_1,\dots,\lambda_n)$, i.e.
chains from $\mathcal{S}_\ell^{-\Sigma}(z_1,\dots,z_n;\lambda_1,\dots,\lambda_n)$
annihilated by the boundary operator.
\end{prop} 

\noindent\underline{\bf Proof.}
At first, let us consider the case of one point $z$, i.e. we will apply the boundary operator to $G_\ell^\ell(z,\lambda)$.

It is clear that 
$\partial  G_\ell^\ell(z,\lambda)=\{pt\}\times\alpha_\ell(q)G_{\ell-1}^{\ell-1}(z,\lambda)$,
where $\{pt\}$ stands for $\infty$-point and $\alpha_\ell(q)$ is some constant.

Analyzing carefully the $q$-factors, we find that
\begin{equation}
G_\ell^\ell(z,\lambda)=\{pt\}\times(\widetilde{E}G_\ell^\ell),
\end{equation}
where
\begin{equation}
\widetilde{E}G_\ell^\ell=\sum_{k=0}^{\ell-1}q^{2(n-k-1)}\Big(1-q^{2\lambda+4k}\Big)G_{\ell-1}^{\ell-1}.
\end{equation}
Each term in the sum above correspond to the application of the boundary operator to each of the $\ell$ loops, while 
$q$-factors appear from different values of $\psi_z$ on the boundaries. From here one can deduce that 
\begin{eqnarray}
\widetilde{E}G_\ell^\ell&=&\varphi_{\vec{z}}\Big(\sum_{k=0}^{\ell-1}q^{2(n-k-1)}(1-q^{-2\lambda+4k})F^{\ell-1}v_\lambda\Big)\nonumber\\
&=&\varphi_{\vec{z}}\big((q-q^{-1})q^HEF^\ell v_\lambda\big).
\end{eqnarray}
Hence, the action of $\widetilde{E}$ and, therefore, of the boundary operator reproduces the action of $E'=(q-q^{-1})q^HE$. 

It is easy to see that 
$\Delta(E')=q^H\otimes E'+E'\otimes 1$. Now let us apply the boundary operator to 
$G_{k_1,k_2}^\ell(z_1,z_2;\lambda_1,\lambda_2)$. We have
\begin{eqnarray}
&&\partial G_{k_1,k_2}^\ell(z_1,z_2;\lambda_1,\lambda_2)=\nonumber\\
&&\{pt\}\times\Big(G_{k_1-1,k_2}^{\ell-1}\alpha_{k_1}(q)+q^{-\lambda_1-2k_1}G_{k_1,k_2-1}^{\ell-1}\alpha_{k_2}(q)\Big)=\nonumber\\
&&\{pt\}\times\varphi_{\vec{z}}\big(\Delta(E')F^{k_1}v_{\lambda_1}\otimes F^{k_2}v_{\lambda_2}\big).
\end{eqnarray}
One can continue it to the case of $n$ points obtaining that the action of 
the boundary operator is equivalent to the action of $(q-q^{-1})q^HE$ on the product of Verma modules.
Therefore, we have one-to-one map between singular vectors in the product of Verma modules and cycles in
$\oplus_\ell\mathcal{S}_\ell^{-\Sigma}(z_1,\dots,z_n;\lambda_1,\dots,\lambda_n)$. 
\hfill$\blacksquare$\bigskip

\noindent In \cite{varchenko}, Varchenko found the homological description of singular vectors in the tensor product of Verma modules. Let us recall the 
notation (\ref{homology}) of subsection 4.1.
\begin{theorem}\label{thvarch} 
There is a natural isomorphism between 
$H^{-\Sigma}_\ell(z_1,z_2...,z_n;$ $\lambda_1,$ $\lambda_2$ $...\lambda_n)$ and singular vectors on the level 
$\lambda-2\ell$ in the tensor product of n Verma modules, i.e. 
$Sing(M_{\lambda_1}\otimes ...\otimes M_{\lambda_n})[\lambda-2\ell]$, where $\lambda=\lambda_1+...+\lambda_n$. 
\end{theorem}
\noindent According to the statement of Theorem \ref{thvarch} and Proposition \ref{cycsing} we can identify $Z^{-\Sigma}_{\ell}(\vec{\lambda},\vec{z})$ with 
$H^{-\Sigma}_{\ell}(\vec{\lambda},\vec{z})$. Then, restricting the result of Proposition \ref{rmatrix} to the subspace of singular vectors and its geometric counterpart,  we rederive the result of Varchenko \cite{varchenko}. 

\begin{theorem}
The following diagram commutes:
 \begin{eqnarray}
\xymatrix{
{\begin{array}{c}
Sing_{\lambda-2\ell}(M_{\lambda_1}\otimes\dots\otimes M_{\lambda_i}\otimes\\
M_{\lambda_{i+1}}\otimes...\otimes M_{\lambda_n})
\end{array}}
\ar[d]^{\check{R}_{i,i+1}} \ar[r]^{i} &
H^{-\Sigma}_\ell(z_1,\dots,z_n;\lambda_1,\dots,\lambda_n) \ar[d]^{\mathscr{A}_{i,i+1}}\\
{\begin{array}{c}
Sing_{\lambda-2\ell}(M_{\lambda_1}\otimes\dots\otimes M_{\lambda_i+1}\otimes\\
M_{\lambda_{i+1}}\otimes ...\otimes M_{\lambda_n})
\end{array}}
\ar[r]^{i} &
H^{-\Sigma}_{\ell}(z_1,\dots,z_n;\lambda_1,\dots,\lambda_n).
}
\end{eqnarray}
Here $i$ is an isomorphism and $\mathscr{A}_{i,i+1}$ is the monodromy operator along the paths:
\begin{eqnarray}
z_i(t)&=&\frac{1}{2}\big((z_i+z_{i+1})+((z_i-z_{i+1})e^{\pi it}\big),\nonumber\\
z_{i+1}(t)&=&\frac{1}{2}\big((z_i+z_{i+1})+((z_{i+1}-z_i)e^{\pi it}\big).
\end{eqnarray}
\end{theorem}
In the case of $\ell=1$, it is easy to see that one of the singular vectors is given by Pochhammer loop $P^1$, which can be easily decomposed in terms of chains from $\mathcal{S}_\ell^{-\Sigma}(z_1,z_2;\lambda_1,\lambda_2)$ (see Fig. 5).
\begin{figure}[hbt] 
\centering           
\includegraphics[width=0.8\textwidth]{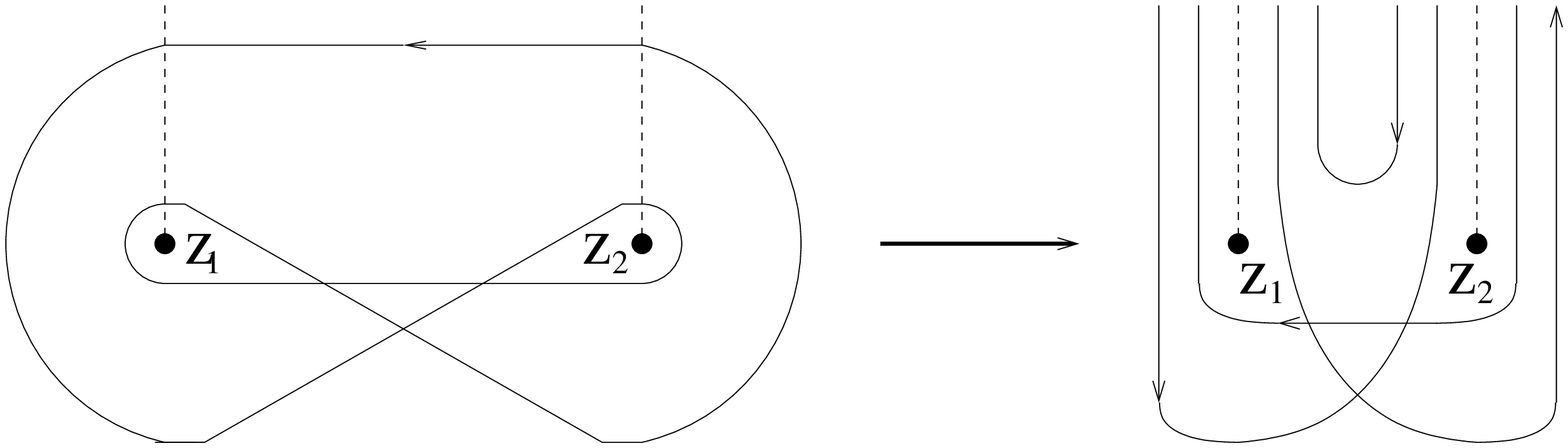}
\caption{Pochhammer loop $P^1$ decomposed}                
\end{figure}     

We remind that in Sec. 2 we identified the space of  $Sing_{\nu}(M_{\lambda_1}\otimes M_{\lambda_2})$ with $Hom(M_{\nu}, M_{\lambda_1}\otimes M_{\lambda_2})$.  Let's consider the cycle from $H_s^{-\Sigma}(z_1,z_2;\lambda_1,\lambda_2)$, corresponding to an element of $Sing_{\lambda_1+\lambda_2-2s}(M_{\lambda_1}\otimes M_{\lambda_2})$. We will denote it by $P_{\lambda_1\lambda_2}^s$ and draw it schematically as the product of $s$ Pochhammer cycles (in case if it is nonzero), motivated by analogy with $s=1$ case. One can notice that there is exactly one $P_{\lambda_1\lambda_2}^s$, i.e. singular vector on the level $\lambda_1+\lambda_2-2s$ if $\lambda_1,\lambda_2\ge 0$ and $0\le s\le \lambda_1+\lambda_2-|\lambda_1-\lambda_2|$.

Now one can construct an intertwiner, using this geometric language and the equivalence between 
the $Hom$-spaces and the spaces of singular vectors. 
\begin{figure}[hbt] 
\centering           
\includegraphics[width=1.0\textwidth]{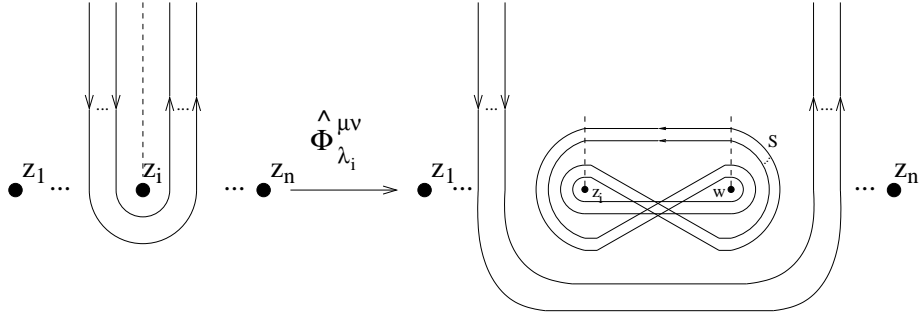}
\caption{Intertwiner's action}                
\end{figure}     
In order to do this, one should consider the map 
$\hat{\Phi}_{\lambda_i}^{\mu\nu}(z,w)$, constructed as it is shown on Fig. 6. 
Namely, this is a map from  
$\mathcal{S}^{-\Sigma}_\ell(z_1,\dots,z_n;\lambda_1,\dots,\lambda_n)$ to 
$\mathcal{S}^{-\Sigma}_{\ell+s}(z_1,\dots,z_i,w,z_{i+1},\dots,z_n;
\lambda_1,\dots,\lambda_{i-1},\mu,\nu,\lambda_{i+1},\dots,\lambda_n)$. It is constructed by means of the 
puncturing of additional point $w$ inside the Gomez-Sierra contours surrounding $z_i$ such that 
$\mathrm{Re}\ z_i<\mathrm{Re}\ w<\mathrm{Re}\ z_{i+1}$ and the insertion of the cycle $P_{\mu\nu}^s$ winding around them. This operator is precisely a geometric version of an intertwiner $\Phi_{\lambda_i}^{\mu\nu}$, since $P_{\mu\nu}^s$ corresponds to a singular vector. Namely, the following Proposition holds.

\begin{prop}\label{intertwiner}
Let $\mathrm{Re}\ z_1<\dots< \mathrm{Re}\ z_n$. Then the following diagram is commutative:
\begin{eqnarray}
\xymatrix{
{\begin{array}{c}
M_{\lambda_1}\otimes\dots\otimes M_{\lambda_i}\otimes\dots\otimes\\
 M_{\lambda_n}[\lambda-2\ell]
\end{array}}
\ar[r]^i \ar[d]^{\Phi_{\lambda_i}^{\mu\nu}}
&
\mathcal{S}^{-\Sigma}_\ell(z_1,\dots,z_n;\lambda_1,\dots,\lambda_n)
\ar[d]^{\hat{\Phi}_{\lambda_i}^{\mu\nu}(z_i,w)}\\
{\begin{array}{c}
M_{\lambda_1}\otimes\dots\otimes M_\mu\otimes M_\nu\otimes\dots\otimes\\
 M_{\lambda_n}[\lambda-2\ell-2s]
\end{array}}
\ar[r]^i
&
{\begin{array}{c}
\mathcal{S}^{-\Sigma}_{\ell+s}(z_1,\dots,z_i,w,z_{i+1},\dots,z_n;\\
\lambda_1,\dots,
\lambda_{i-1},\mu,\nu,\lambda_{i+1},\dots,\lambda_n),
\end{array}}
}
\end{eqnarray}
where $\Phi_{\lambda_i}^{\mu\nu}$ is the element of $Hom(M_{\lambda_i},M_\mu\otimes M_\nu)$ corresponding to the singular vector  
from $Sing_{\lambda_i}(M_{\mu}\otimes M_{\nu})$ geometrically represented by $P^s_{\lambda_1\lambda_2}$, 
and  $\mathrm{Re}\ z_i<\mathrm{Re}\ w<\mathrm{Re}\ z_{i+1}$.
\end{prop}

\noindent Therefore, for generic values of weights we have the following Proposition which gives the bilinear relations between "geometric" intertwiners, defined above. 
\begin{prop}\label{genint}
Let $\mathrm{Re}\ z_1< \mathrm{Re}\ z_2<\mathrm{Re}\ z_3$ and $\lambda_i$ $(i=0,1,2,3)$ be generic. Then we have
\begin{eqnarray}
&&\mathscr{A}_{2,3}\big(\hat{\Phi}_\rho^{\lambda_1\lambda_2}(z_1,z_2)\hat{\Phi}_{\lambda_0}^{\rho\lambda_3}
(z_1,z_3)\big)=\nonumber\\
&&=\sum_{\xi}B^M_{\rho\xi}
\left[\begin{array}{cc}
\lambda_0 & \lambda_1\\
\lambda_2 & \lambda_3
\end{array}
\right]
\hat{\Phi}_\xi^{\lambda_1\lambda_3}(z_1,z_3)\hat{\Phi}_{\lambda_0}^{\xi\lambda_2}(z_1,z_2),
\end{eqnarray}
where the intertwiners $\hat{\Phi}^{\mu}_{\nu\lambda}$ in the expression above act on the space $\mathcal{S}^{-\Sigma}_\ell(\vec{z};\vec{\lambda})$.
\end{prop}
\noindent\underline{\bf Proof.} This statement is a consequence of Propositions \ref{intertwiner0} and \ref{rmatrix}. 
We already know, that the action of the monodromy operator is equivalent to the action of the R-matrix on the product of modules. Therefore, 
$\mathscr{A}_{2,3}\big(\hat{\Phi}_\rho^{\lambda_1\lambda_2}(z_1,z_2)$ $\hat{\Phi}_{\lambda_0}^{\rho\lambda_3}(z_1,z_3)\big)$ corresponds to the expression 
$(1\otimes PR)\Phi_\rho^{\lambda_1\lambda_2}\Phi_{\lambda_0}^{\rho\lambda_3}$. On the other hand, we know that this is equal (see (\ref{int1})) to 
$\sum_{\xi}B^M_{\rho\xi}\left[\begin{array}{cc}
\lambda_0 & \lambda_1\\
\lambda_2 & \lambda_3
\end{array}
\right]
\Phi_\xi^{\lambda_1\lambda_3}\Phi_{\lambda_0}^{\xi\lambda_2}.$ Thus the Proposition is proven.
\hfill$\blacksquare$\bigskip

\noindent {\bf 4.4.  Irreducible representations of $U_q(sl(2))$ via local systems.} In this section, we considered relations between tensor products of Verma modules and local systems on configuration space. There is also a possibility to include the irreducible modules of $U_q(sl(2))$ in such a correspondence. In order to do that one should understand the factorization map $M_{\lambda}\to M_{\lambda}/M_{-\lambda-2}$ (for $\lambda \ge 0$) in the geometric context. This can be established by considering the relative cycles (w.r.t. the set of $x_i=z_j$) instead of usual ones.  

Namely, we denote by $\tilde {G}_{r_1,r_2\dots r_n}^\ell [z_{1},...,z_{n}, \lambda_1,...,\lambda_n]$ the same geometric object as before, but consider it as a relative cycle not only with respect to the hyperplanes corresponding to the reference point, but also 
with respect to the set $x_i=z_{p}$, where $i=1,..,\ell$ and $p=1,...,n$. Let us denote the antisymmetric (w.r.t. the action of the group $\Sigma$) part of the space of such cycles factorized by homologically trivial ones  
as  $\tilde{\mathcal{S}}^{-\Sigma}_\ell [z_1,\dots,z_n;\lambda_1,\dots,\lambda_n]$. 
One can give an explicit geometric picture of transfer from the representatives corresponding to absolute cycles (associated with basic  contours) to relative  cycles. 
It can be achieved by means of the shrinking of the appropriate contours in $G_{r_1,r_2\dots r_n}^\ell (z_{1},...,z_{n}, \lambda_1,...,\lambda_n)$  and the moving all the resulting rays from $z_{k}$ to $\infty$ on the left-hand side of the cut, as it is shown in Fig. 7.
Counting the $q$-powers at the points $z_{1}, ...., z_{n}$, we obtain that the elements  
represented by $\tilde {G}_{r_1,r_2\dots r_n}^\ell [z_{1},...,z_{n}]\in\tilde{\mathcal{S}}^{-\Sigma}_\ell[z_1,\dots,z_n;\lambda_1,\dots,$ $\lambda_n]$, will be annihilated if $r_{i}\ge \lambda_{i}$ and $\lambda_{i}\in \mathbb{Z}_+$ $(i=1,\dots, n)$. In the example of one point, as shown on Fig. 7, the resulting factor, corresponding to the right hand side of the picture, is: $\prod_{k}(1-q^{2(\lambda-k)})$ (each multiplier is a contribution from a different basic contour), which vanishes, when $k$ reaches $\lambda$. In such a way we have a natural map between 
$\oplus_\ell \tilde{\mathcal{S}}_\ell^{-\Sigma}[z_1,\dots, z_n;\lambda_1,\dots,\lambda_n]$ and 
$V_{\lambda_1}\otimes\dots\otimes V_{\lambda_n}$. Namely, it is the map:
\begin{eqnarray}
\xymatrix{V_{\lambda_1}\otimes\dots\otimes V_{\lambda_n}[\lambda-2\ell]
\ar[r]^{\tilde{\varphi}_{\vec{z}}} &
\tilde{\mathcal{S}}^{-\Sigma}_\ell [z_1,\dots,z_n;\lambda_1,\dots,\lambda_\ell],\\
}
\end{eqnarray}
such that for $k_1+\dots+k_n=\ell$
\begin{equation}
\tilde{\varphi}_{\vec{z}}:\ F^{k_1}v_{\lambda_1}\otimes\dots\otimes F^{k_n}v_{\lambda_n}
\longmapsto
\tilde{G}^{\ell}_{k_1,\dots,k_n}[z_1,\dots,z_n;\lambda_1,\dots,\lambda_\ell].
\end{equation}

\begin{figure}[hbt] 
\centering           
\includegraphics[width=0.5\textwidth]{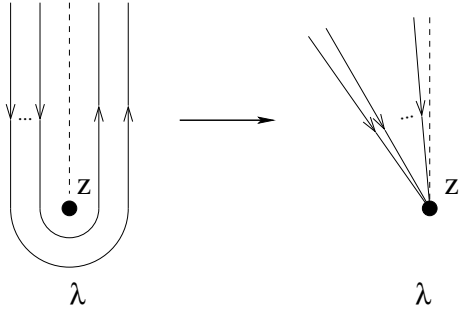}
\caption{Shrinking of a cycle}                
\end{figure}     
This allows us to formulate a statement. 
\begin{prop}
Let $\lambda_{i}\in \mathbb{Z}_+$ $(i=1,\dots, n)$. Then the following diagram is commutative:
\begin{eqnarray}
\xymatrix{
M_{\lambda_1}\otimes\dots\otimes M_{\lambda_n}[\lambda-2\ell]      
\ar[d]^{p_{1...n}} \ar[r]^{ \qquad\qquad\varphi_{\vec{z}}}&
\mathcal{S}^{-\Sigma}_\ell(\vec{z},\vec{\lambda})
   \ar[d]^{f_{1...n}}\\
{V_{\lambda_{1}}\otimes ... \otimes V_{\lambda_{n}}[\lambda-2\ell]}
         \ar[r]^{\qquad\qquad\tilde{\varphi}_{\vec{z}}} & 
         \tilde{\mathcal{S}}^{-\Sigma}_\ell[\vec{z};\vec{\lambda}],}
\end{eqnarray} 
where the map $p_{1...n}$ is just the composition of projections on the corresponding irreducible modules and the map $f_{1...n}$ is given by the composition of shrinking maps as demonstrated on Fig. 7.
\end{prop}
\noindent Similar construction of tensor product of finite-dimensional representations was used in \cite{fkv} in the framework of \cite{varchenko}. 

Our last task will be to give the geometric meaning to the relations between intertwiners, corresponding to finite-dimensional representations.

\begin{prop}\label{fdgeom}
Let $\mathrm{Re}\ z_1< \mathrm{Re}\ z_2<\mathrm{Re}\ z_3$ and $\lambda_i\in\mathbb{Z}_+$ $(i=0,1,2,3)$. Then we have
\begin{eqnarray}\label{mono}
&&\mathscr{A}_{2,3}\big(\hat{\Phi}_\rho^{\lambda_1\lambda_2}(z_1,z_2)\hat{\Phi}_{\lambda_0}^{\rho\lambda_3}
(z_1,z_3)\big)=\nonumber\\
&&=\sum_{\xi}B^V_{\rho\xi}
\left[\begin{array}{cc}
\lambda_0 & \lambda_1\\
\lambda_2 & \lambda_3
\end{array}
\right]
\hat{\Phi}_\xi^{\lambda_1\lambda_3}(z_1,z_3)\hat{\Phi}_{\lambda_0}^{\xi\lambda_2}(z_1,z_2),
\end{eqnarray}
where $\rho,\xi\in \mathbb{Z}_+$ and the intertwiners $\hat{\Phi}^{\mu}_{\nu\lambda}$ in the expression above act on the space $\tilde{\mathcal{S}}^{-\Sigma}_\ell[\vec{z};\vec{\lambda}]$.
\end{prop}
Propositions \ref{genint} and \ref{fdgeom} are the geometric versions of Propositions \ref{intertwiner0} and \ref{intertwiner1}. 
Formula (\ref{mono}) is very important for deriving the relation between the intertwiners of Fock space modules, which we will study in the next section.

\section{The construction of intertwiners between Fock spaces and braided VOA on \\
${\mathbb{F}_\varkappa=\bigoplus_{\lambda\ge 0}(V_{\Delta(\lambda),\varkappa}\otimes V_\lambda)}$}
\setcounter{equation}{0}
\noindent{\bf {5.1. Intertwiners between Fock spaces.}} Let's consider the multivalued function $\Psi_{\vec{z}}(x_1,\dots,x_s)$ (\ref{psi}). As we already know, this function determines a local system on a configuration space.
One can consider differential forms on $\mathbb{C}^s$ of the following kind:
\begin{equation}\label{dform}
\Psi_{z_1,\dots,z_n}(x_1,\dots,x_s)A(z_1,\dots,z_n;x_1,\dots,x_s)\ud x^1\wedge\dots\wedge\ud x^s
\end{equation}
and integrate them over cycles in $H_s^{-\Sigma}(\mathbb{C}^s\backslash \mathcal{H},\mathscr{S})$, if
$A(z_1,\dots,z_n;x_1,\dots,x_s)$ is symmetric in $x_i$.
Now we reduce the number of $z$-variables to two of them, namely, we choose them to be $z,0$ and 
consider the following expression for $A$:
\begin{equation}\label{fint}
A^{\lambda,v,v^*}(z,0;x_1,\dots,x_s)=<v^*,:\mathbb{X}(\lambda,z)\mathbb{X}_s^+(x_1)\dots\mathbb{X}_s^+(x_s):v>,
\end{equation}
where $v\in F_{\mu,\varkappa}\ $ and $v^*\in F^*_{\lambda+\mu-2s,\varkappa}$.
If we integrate the corresponding differential form over a cycle $P_{\lambda\mu}^s\in H_s^{-\Sigma}(\mathbb{C}^s\backslash \mathcal{H},\mathscr{S})$, which we constructed in the previous section, one can consider the resulting expression
as a matrix of some operator
\begin{equation}
\varphi(\cdot,z): F_{\mu,\varkappa}\rightarrow F_{\lambda+\mu-2s,\varkappa}.
\end{equation}
Moreover, the following statement holds \cite{ff2}.

\begin{prop}
Let $\lambda,\mu, \nu\in\mathbb{Z}$ such that $\nu\le \lambda+\mu$. Then there exists an intertwining operator
\begin{equation}\label{intvir}
\Phi_{\lambda\mu}^\nu(z):F_{\lambda,\varkappa}\otimes F_{\mu,\varkappa}\rightarrow 
F_{\nu,\varkappa}[[z,z^{-1}]]z^{\Delta(\nu)-\Delta(\mu)-\Delta(\lambda)},
\end{equation}
i.e. the operator, such that
\begin{equation}\label{propint}
L_n\cdot\Phi_{\lambda\mu}^\nu(z)=\Phi_{\lambda\mu}^\nu(z)\Delta_{z,0}(L_n),
\end{equation}
where
\begin{equation}
\Delta_{z,0}(L_n)=\oint_z\frac{\ud\xi}{2\pi i}\xi^{n+1}\Big(\sum_m(\xi-z)^{-m-2}L_m\Big)\otimes 1+1\otimes L_n.
\end{equation}
In particular case when the first argument is the highest weight vector  
${\bf 1}_{\lambda}\in F_{\lambda,\varkappa}$, the explicit expression for the matrix elements of an intertwiner are given by the following formula:
\begin{eqnarray}
&&\langle v^*,\Phi_{\lambda\mu}^\nu(z)({\bf{1}}_{\lambda}\otimes v)\rangle=
\int_{P_{\lambda\mu}^s}\Psi_{0,z}(x_1,\dots,x_s)
\nonumber\\
&&\langle v^*,:\mathbb{X}(\lambda,z)\mathbb{X}_s^+(x_1)\dots\mathbb{X}_s^+(x_s)
:v\rangle\ud x^1\wedge\dots\wedge\ud x^s,
\end{eqnarray}
where  $v\in F_{\mu,\varkappa}\ $, and $v^*\in F^*_{\lambda+\mu-2s,\varkappa}$
($s=\frac{\lambda+\mu-\nu}{2}$).
\end{prop}

\noindent\underline{\bf Proof.}
One can construct an intertwining operator by means of the following procedure. Let's consider the  correlator
\begin{eqnarray}
\langle v^*,Y(u,z)\mathbb{X}_s^+(x_1)\dots\mathbb{X}_s^+(x_s)v\rangle,
\end{eqnarray}
where $u\in F_{\lambda,\varkappa}$ and $z,x_1,,,x_s\in \mathbb{R}$ such that $z>x_1>...>x_s$. 
This expression can be rewritten in the form 
\begin{eqnarray}
\Psi_{0,z}(x_1,...,x_s)f^u_{v,v^*}(z,x_1,...,x_s),
\end{eqnarray}
where $f^u_{v,v^*}(z,x_1,...,x_s)$ is a rational function of $z,x_1,\dots, x_s$. One can make the analytic 
continuation of this multivalued function to the complex domain, using the branches of 
$\Psi_{0,z}(x_1,...,x_s)$. Therefore, one can define the matrix elements of the intertwiner by 
 \begin{eqnarray}\label{explint}
&&\langle v^*,\Phi_{\lambda\mu}^\nu(z)(u\otimes v)\rangle=\nonumber\\
&&\int_{P_{\lambda\mu}^s}\Psi_{0,z}(x_1,\dots,x_s)f^u_{v,v^*}(z,x_1,...,x_s)\ud x^1\wedge\dots\wedge\ud x^s.
\end{eqnarray}
In order to check the property (\ref{propint}), one needs to consider the expression
\begin{eqnarray}
\langle v^*, \int_C \frac{\ud \xi }{2\pi i}\xi^{n+1}L(\xi)Y(u,z)\mathbb{X}_s^+(x_1)\dots\mathbb{X}_s^+(x_s)v\rangle,
\end{eqnarray}
where the contour $C$ is a circle including point $0$ and $z$. The formula above can be reexpressed in the following way:
\begin{eqnarray}
&&\langle v^*, \oint_z \frac{\ud \xi }{2\pi i}\xi^{n+1}(\xi-z)^{-n-2}Y(L_n\cdot u,z)\mathbb{X}_s^+(x_1)\dots\mathbb{X}_s^+(x_s)v +\nonumber\\
&&Y(u,z)\mathbb{X}_s^+(x_1)\dots\mathbb{X}_s^+(x_s)L_n\cdot v+\nonumber\\
&&\sum^s_{i=1}Y(u,z)\mathbb{X}_s^+(x_1)\dots \partial_{x_i}(x_i^{n+1}\mathbb{X}_s^+(x_i)) \dots\mathbb{X}_s^+(x_s)\cdot v\rangle.
\end{eqnarray}
Since in the formula (\ref{explint}) we are integrating over closed cycle, terms with derivatives of screening operators disappear due to Stokes theorem and we arrive to the formula (\ref{propint}).

One can obtain the coefficient $z^{\Delta(\nu)-\Delta(\mu)-\Delta(\lambda)}$
by simple change of variables. Namely, replacing $x_i\rightarrow\frac{x_i}{z}$
we get an expression for the nonsingle-valued multiplier from 
$<v^*,\Phi_{\lambda\mu}^\nu(z)(u\otimes v)>$:
\begin{equation}
(z_1-z_2)^{\Delta_{\lambda\mu}^\nu}, \qquad \Delta_{\lambda\mu}^\nu=s+\frac{\lambda\mu}{2\varkappa}+
\frac{s(s-1)}{\varkappa}-\frac{\lambda+\mu}{2\varkappa}.
\end{equation} 
One can easily see that $\Delta_{\lambda\mu}^\nu=\Delta(\nu)-\Delta(\lambda)-\Delta(\mu)$.
\hfill$\blacksquare$\bigskip

As in the case of $U_q(sl(2))$ we will be interested in the composition of intertwining operators. 
The composition of intertwiners $\Phi_{\lambda_3\rho}^{\lambda_0}(z)\Phi_{\lambda_2\lambda_1}^\rho(w)$
is understood as an operator, acting on $F_{\lambda_3}\otimes F_{\lambda_2}\otimes F_{\lambda_1}$.

We remind that when we studied local systems, we have introduced the operators
$\hat\Phi_\lambda^{\mu \nu}(z_1,z_2)$ geometrically representing intertwining operators.
They obey the braiding relation, repeating the one for the 
finite-dimensional representations of $U_q(sl(2))$. Namely, for 
$0<\mathrm{Re}\ z_1<\mathrm{Re}\ z_2$
\begin{eqnarray}\label{intrel}
&&\mathscr{A}_{z_1,z_2}\big(\hat\Phi_\rho^{\lambda_1\lambda_2}(0,z_1)
\hat\Phi_{\lambda_0}^{\rho\lambda_3}(0,z_2)\big)=\nonumber\\
&&=\sum_{\xi}B^V_{\rho\xi}
\left[\begin{array}{cc}
\lambda_0 &\lambda_1 \\
\lambda_2 &\lambda_3 
\end{array}
\right]
\hat\Phi_\xi^{\lambda_1\lambda_3}(0,z_2)
\hat\Phi_{\lambda_0}^{\xi\lambda_2}(0,z_1),
\end{eqnarray}
where the expression above acts on $\oplus_\ell \tilde{S}_\ell^{-\Sigma}[\lambda_1,\lambda_2,\lambda_3;0,z_1,z_2]$ and 
$\lambda_i\in \mathbb{Z_+}$, $(i=0,1,2,3)$. 
If one integrates suitable expressions, like the integrand of (\ref{explint}) over the cycles from the expression above, one arrives to the following statement. 
\begin{prop}
Let $z_1,z_2\in \mathbb{R}$, such that $0<z_1<z_2$ and $\lambda_i\ge 0$ $(i=0,1,2,3)$, 
$\varkappa>{\rm max}{\lambda_i}$. Then the following relation holds:
\begin{eqnarray}\label{fockint}
&&\mathscr{A}_{z_1,z_2}\big(\Phi_{\lambda_3\rho}^{\lambda_0}(z_2)\Phi_{\lambda_2\lambda_1}^\rho(z_1)\big)(P\otimes 1)=\nonumber\\
&&\sum_{\xi}B^V_{\rho\xi}
\left[\begin{array}{cc}
\lambda_0 &\lambda_1 \\
\lambda_2 &\lambda_3 
\end{array}
\right]
\Phi_{\lambda_2\xi}^{\lambda_0}(z_1)\Phi_{\lambda_3\lambda_1}^\xi(z_2), 
\end{eqnarray} 
where $P$ is the interchange operator, namely $P(v_1\otimes v_2)=v_2\otimes v_1$.
\end{prop}

Using the explicit values for the "screened" correlators in Coulomb gas formalism studied in the physics literature (see e.g. \cite{df1}, \cite{df2}), one can show that the relation (\ref{fockint}) can be analytically continued in $\varkappa$ to any value in $\mathbb{R}\backslash\mathbb{Q}$.\\

\noindent{\bf{5.2. Braided VOA on $\mathbb{F}_\varkappa =$$\bigoplus_\lambda(V_{\Delta(\lambda),\varkappa}$$\otimes V_\lambda)$}}. 
 Let's consider the following space: $\mathbb{F}_\varkappa =$$\bigoplus_{\lambda\in \mathbb{Z}_+}(V_{\Delta(\lambda),\varkappa}$$\otimes V_\lambda)$.  Below we will show that there exists a structure of a braided VOA on this space.

First of all, we define the following map:
\begin{equation}\label{bvoa}
Y:v\otimes a\rightarrow Y(v\otimes a,z)=
\sum_{\nu,\mu} \Phi^{\nu}_{\lambda\mu}(z)(v\otimes \cdot)\otimes \phi^{\nu}_{\mu\lambda}(\cdot \otimes a).
\end{equation}
Here $v\in F_{\lambda,\varkappa}$ and $a\in V_\lambda$ for  some $\lambda\in\mathbb{Z}$. 
As a consequence of the Proposition \ref{qminus} we have the following statement.
\begin{lemma}
$[Q^-,\Phi^{\nu}_{\mu\lambda}(z)(v\otimes \cdot)]=0$ if 
$v\in V_{\Delta(\lambda),\varkappa}\subset F_{\lambda,\varkappa}$.
Hence, $[Q^-\otimes 1,Y(v\otimes a,z)]=0$ if $v\in V_{\Delta(\lambda),\varkappa}$.
\end{lemma}
\noindent Therefore, $Y$ acts as follows:
\begin{equation}
Y: \mathbb{F}_{\varkappa}\to 
End(\mathbb{F}_{\varkappa})\{z\}.
\end{equation}
Let $0<z_1<z_2$, $z_1,z_2\in \mathbb{R}$, $v_i\in V_{\Delta{(\lambda_i)},\varkappa}$ $(i=1,2,3)$, 
$a_i\in V_{\lambda_i}$ $(i=1,2,3)$. Then 
\begin{eqnarray}
&&\mathscr{A}_{z_2,z_1}\big(Y(v_1\otimes a_1,z_2)Y(v_2\otimes a_2,z_1)\big)(v_3\otimes a_3)=\nonumber\\
&&\mathscr{A}_{z_2,z_1}\Big(\sum_{\lambda_1,\lambda_2,\nu,\rho}
\Phi_{\lambda_1\rho}^\nu(z_2)\Phi_{\lambda_2\lambda_3}^\rho(z_1)\otimes
\phi_{\rho\lambda_1}^\nu\phi_{\lambda_3\lambda_2}^\rho\Big)\cdot\nonumber\\
&&(v_1\otimes v_2\otimes v_3)\otimes(a_3\otimes a_2\otimes a_1)=\nonumber\\
&&\Big(\sum_{\lambda_1,\lambda_2,\nu,\rho,\xi}\Phi_{\lambda_2\xi}^\nu(z_1)\Phi_{\lambda_1\lambda_3}^\xi(z_2)
B^V_{\rho\xi}
\left[\begin{array}{cc}
\nu &\lambda_3 \\
\lambda_2 &\lambda_1 
\end{array}
\right] \otimes 
\phi_{\rho\lambda_1}^\nu\phi_{\lambda_3\lambda_2}^\rho\Big)\cdot\nonumber\\
&&(v_2\otimes v_1\otimes v_3)(a_3\otimes a_2\otimes a_1)=\nonumber\\
&&\Big(\sum_{\lambda_1,\lambda_2,\nu,\xi}
\Phi_{\lambda_2\xi}^{\nu}(z_1)\Phi_{\lambda_2\lambda_3}^\xi(z_2)\otimes
\phi_{\xi\lambda_2}^{\nu}\phi_{\lambda_3\lambda_1}^\xi\Big)\cdot\nonumber\\
&&(v_2\otimes v_1\otimes v_3)\otimes(a_3\otimes \sum_ir_i^{(2)}a_1\otimes r_i^{(1)}a_2)=\nonumber\\
&&\sum_i Y(v_2\otimes r_i^{(1)}a_2,z_1)Y(v_1\otimes r_i^{(2)}a_1,z_1)(v_3\otimes a_3).
\end{eqnarray}
Hence, we have proved the following Proposition.

\begin{prop}\label{commu}
Map $Y$ satisfies commutativity condition, namely let $z,w\in \mathbb{R}$, such that 
$0<z<w$, then
\begin{eqnarray}
\mathscr{A}_{z,w}\big(Y(v_1\otimes a_1,w)Y(v_2\otimes a_2,z)\big)=\nonumber\\
\sum_i Y(v_2\otimes r_i^{(1)}a_2,z)Y(v_1\otimes r_i^{(2)}a_1,w),
\end{eqnarray} 
where $R=\sum_i  r_i^{(1)}\otimes r_i^{(2)}$ is the universal R-matrix for $U_q(sl(2))$.
\end{prop}

The next proposition shows that  the Virasoro action is compatible with correspondence $Y$.

\begin{prop}
\hfill
\begin{itemize}
\item[(i)]There is a natural action of Virasoro algebra on $ \tilde {\mathbb{F}}_{\kappa}=
\oplus_{\lambda\in \mathbb{Z}_{+}} V_{\Delta{(\lambda)},\varkappa}\otimes V_\lambda^c$, namely 
$L_n\cdot(v\otimes a)=(L_n\cdot v)\otimes a$, 
where $v\in F_{\lambda,\varkappa}, a\in M_\lambda^c$ for some $\lambda$.

Moreover, there exists an element $\tilde{\omega}\in V_{\Delta(0),\varkappa}\otimes V_0$,
namely $\tilde{\omega}=\omega\otimes v_0$ ($v_0$ is the only basis element in $V_0$),
such that $Y(\tilde{\omega},z)=L(z)=\sum_n L_nz^{-n-2}$ is the Virasoro element.

\item[(ii)] $[L_{-1},Y(v\otimes a,z)]=\partial_z Y(v\otimes a,z)$.
\end{itemize}
\end{prop}

\noindent\underline{\bf Proof.}
(i) is obvious.
(ii) follows from the explicit definition of the intertwining operator $\Phi_{\lambda\mu}^{\nu}(z)$.
\hfill$\blacksquare$\bigskip

\noindent Finally, we show that the associativity condition holds (see e.g \cite{FHL}). First, we will prove the following Lemma.

\begin{lemma} \label{der}
Map Y satisfies creation property, namely:
\begin{eqnarray*}
&& (i)\ Y(v\otimes a,z)\mathbf{1}=e^{zL_{-1}}(v\otimes a),\\
&& (ii)\ e^{z_2L_{-1}}Y(v\otimes a,z_1) e^{-z_2L_{-1}}=Y(v\otimes a,z_1+z_2).
\end{eqnarray*}
\end{lemma} 

\noindent\underline{\bf Proof.}
At first, we prove (ii).
Let $\mathrm{ad}_{L_{-1}}\cdot=[L_{-1},\cdot]$.
Then $\mathrm{ad}_{L_{-1}}Y(v\otimes a,z)=\partial_w^nYv\otimes a,z)$.
Now
\begin{equation}
\sum_{n=0}^\infty\frac{z_2^n}{n!}\partial_{z_1}^nY(v\otimes a,z_1)=
Y(v\otimes a,z_1+z_2).
\end{equation}
Therefore,
\begin{equation}
e^{z_2L_{-1}}Y(v\otimes a,z_1) e^{-z_2L_{-1}}=Y(v\otimes a,z_1+z_2).
\end{equation}
(i) follows as an easy consequence.
\hfill$\blacksquare$\bigskip

\begin{prop}
Let $t,w,z \in \mathbb{R}$, such that $0< t<w<z$. Then
\begin{eqnarray}\label{asso}
&&Y(v_1\otimes a_1,z)Y\big( Y(v_2\otimes a_2,w-t)v_3\otimes a_3,t\big)\mathbf{1}=\nonumber\\
&&Y\big( Y(v_1\otimes a_1,z-w)v_2\otimes a_2,w \big)Y(v_3\otimes a_3,t)\mathbf{1}.
\end{eqnarray}
\end{prop}

\noindent\underline{\bf Proof.} At first, we remind the quasitriangular 
property of the universal R-matrix:
\begin{equation}
(I\otimes \Delta)R=R^{13}R^{12},
\end{equation}
or, in components:
\begin{equation}
r_i^{(1)}\otimes \Delta (r_i^{(2)})=\sum_{i,j}r_i^{(1)}r_j^{(2)}\otimes r_j^{(2)}\otimes r_i^{(2)}.
\end{equation}
We combine it with Lemma \ref{der} and Proposition \ref{commu} and derive (\ref{asso}) as follows.
Let us denote by $\widetilde{\mathscr{A}}_{z,w}$ the inverse of our usual monodromy:
\vspace{3mm}

\begin{eqnarray}
\xymatrix{
 z\ \bullet\ \ar@/_1pc/[r] &
\bullet \   \ar@/_1pc/[l] w
}
\end{eqnarray}

\vspace{3mm}
\noindent Then,
\begin{eqnarray}
&&Y(v_1\otimes a_1,z)Y(v_2\otimes a_2,w)Y(v_3\otimes a_3,t)\mathbf{1}=\nonumber\\
&&\widetilde{\mathscr{A}}_{t,w}\Big(\sum_jY(v_1\otimes a_1,z)
Y(v_3\otimes r_j^{(1)}a_3,t)Y(v_2\otimes r_j^{(2)}a_3,w)\Big)\mathbf{1}=\nonumber\\
&&\widetilde{\mathscr{A}}_{t,z}\widetilde{\mathscr{A}}_{t,w}
\Big(\sum_jY(v_3\otimes r_j^{(1)}r_j^{(2)}a_3,t)Y(v_1\otimes r_j^{(2)}a_1,z)\cdot\nonumber\\
&&Y(v_2\otimes r_j^{(2)}a_2,w)\mathbf{1}\Big)=\nonumber\\
&&\widetilde{\mathscr{A}}_{t,z}\widetilde{\mathscr{A}}_{t,w}
\Big(\sum_{i,j}Y(v_3\otimes r_i^{(1)}r_j^{(1)}a_3,t)Y(v_1\otimes r_i^{(2)}a_1,z)
e^{w L_{-1}}v_2\otimes r_j^{(2)}a_2\Big)=\nonumber\\
&&\widetilde{\mathscr{A}}_{t,z}\widetilde{\mathscr{A}}_{t,w}
\Big(\sum_{i,j}Y(v_3\otimes r_i^{(1)}r_j^{(1)}a_3,t)e^{w L_{-1}}
Y(v_1\otimes r_i^{(2)}a_1,z-w)v_2\otimes r_j^{(2)}a_2\Big)\mathbf{1}=\nonumber\\
&&\widetilde{\mathscr{A}}_{t,z}\widetilde{\mathscr{A}}_{t,w}
\Big(\sum_{i,j}Y(v_3\otimes r_i^{(1)}r_j^{(1)}a_3,t)
Y(Y(v_1\otimes r_i^{(2)}a_1,z-w)v_2\otimes r_j^{(2)}a_2,w)\Big)\mathbf{1}=\nonumber\\
&&=Y\big(Y(v_1\otimes a_1,z-w)v_2\otimes a_2,w\big)Y(v_1\otimes a_3,t)\mathbf{1}.
\end{eqnarray}
In the last line, we used the quasitriangular property.
\hfill$\blacksquare$\bigskip

All these properties allow us to give a general definition of general braided VOA (cf. Proposition \ref{brvoab} and see also 
\cite{dissert}).\\

\noindent{\bf Definition.} {\it Let $\mathbb{V}=\oplus_{\lambda\in I}\mathbb{V}_{\lambda}$ be a direct sum of graded complex 
vector spaces, called sectors: $\mathbb{V}_{\lambda}=\oplus_{n\in \mathbb{Z}_+}\mathbb{V}_{\lambda}[n]$, indexed by some set $I$. Let $\Delta_{\lambda}$, $\lambda\in I$ be complex numbers, which we will call conformal weights of the corresponding sectors.
We say that $\mathbb{V}$ is a braided vertex operator algebra, if there are distinguished elements $0\in I$ such  that $\Delta_0=0$, $\mathbf{1}\in \mathbb{V}_0[0]$, linear maps $D:\mathbb{V}\to\mathbb{V}$, $\mathcal{R}:\mathbb{V}\otimes\mathbb{V}\to \mathbb{V}\otimes\mathbb{V}$ and the linear correspondence 
\begin{eqnarray}
\mathbb{Y}(\cdot,z)\cdot:\mathbb{V}\otimes\mathbb{V}\to \mathbb{V}\{z\},\quad \mathbb{Y}=\sum_{\lambda,\lambda_1,\lambda_2}\mathbb{Y}^{\lambda_1\lambda_2}_{\lambda}(z),
\end{eqnarray}
where 
\begin{eqnarray}
\mathbb{Y}^{\lambda_1\lambda_2}_{\lambda}(z)\in Hom(\mathbb{V}_{\lambda_1}\otimes\mathbb{V}_{\lambda_2},\mathbb{V}_{\lambda})\otimes z^{\Delta_{\lambda}-\Delta_{\lambda_1}-\Delta_{\lambda_2}}
\mathbb{C}[[z,z^{-1}]],
\end{eqnarray}
such that the following properties are satisfied:\\
i)Vacuum property: $\mathbb{Y}(\mathbf{1},z)v=v$, $\mathbb{Y}(v,z)\mathbf{1}|_{z=0}=v$.\\
ii) Complex analyticity: for any $v_i\in \mathbb{V}_{\lambda_i}$, $(i=1,2,3,4)$ the matrix elements 
$\langle v_4^*, \mathbb{Y}(v_3,z_2)\mathbb{Y}(v_2,z_1)v_1 \rangle$ regarded as formal Laurent series in $z_1,z_2$, converge in the domain $|z_2|> |z_1|$ to a complex analytic function $r(z_1,z_2)\in z_1^{h_1}z_2^{h_2}
(z_1-z_2)^{h_3} \mathbb{C}[z_1^{\pm 1},z_2^{\pm 1}, (z_1-z_2)^{-1}]$, where 
$h_1, h_2,h_3\in \mathbb{C}$.\\
iii) Derivation property: $\mathbb{Y}(Dv,z)\mathbf{1}=\frac{d}{dz}\mathbb{Y}(v,z).$\\
iv) Braided commutativity (understood in a weak sense):
\begin{eqnarray}
\mathscr{A}_{z_1,z_2}(\mathbb{Y}(v,z_1)\mathbb{Y}(u,z_2))=\sum_i\mathbb{Y}(u_i,z_2)\mathbb{Y}(v_i,z_1),
\end{eqnarray}
where $\mathcal{R}(u\otimes v)=\sum_i u_i\otimes v_i$.\\
v) There exists an element $\omega\in \mathbb{V}_0$, such that
\begin{equation}
Y(\omega,z)=\sum_{n\in\mathbb{Z}}L_nz^{-n-2}
\end{equation}
and $L_n$ satisfy the relations of Virasoro algebra with $L_{-1}=D$.\\
vi) Associativity (understood in a weak sense):
\begin{eqnarray}
\mathbb{Y}(\mathbb{Y}(u,z_1-z_2)v, z_2)=\mathbb{Y}(u,z_1)\mathbb{Y}(v,z_2).
\end{eqnarray}}\\

\noindent In such a way, we have constructed the operators satisfying all necessary properties 
of VOA algebra on $\mathbb{F}_{\varkappa}$. Therefore, we proved the following statement.
\begin{theorem}
The correspondence 
$Y:\mathbb{F}_\varkappa\rightarrow End(\mathbb{F}_\varkappa)\{z\}$ defined by (\ref{bvoa})
gives a braided VOA structure on $\mathbb{F}_\varkappa$.
\end{theorem}

\section{Identification of the semi-infinite cohomology for $\mathbb{F}_\varkappa\otimes 
\mathbb{F}_{-\varkappa}$} 
\setcounter{equation}{0}
{\bf 6.1. Semi-infinite cohomology: a reminder.} 
In this subsection, we just remind basic facts about semi-infinite cohomology for Virasoro algebra and provide some basic statements which allow to compute explicit formulas for cycles.\\
In the special case of Virasoro algebra, semi-infinite forms can be realized by means of the following super Heisenberg algebra:
\begin{equation}
\{b_n, c_m\}=\delta_{n+m,0}, \quad n,m \in \mathbb{Z}.
\end{equation}
One can construct a Fock module $\Lambda$ in such a way:
\begin{eqnarray}
\Lambda&=&\{b_{-n_1}\dots b_{-n_k}c_{-m_1}\dots c_{-m_\ell} \mathbf{1};\nonumber\\
&& c_k \mathbf{1}=0,\ k\geqslant 2;\quad b_k\mathbf{1}=0, \ k\geqslant -1\}.
\end{eqnarray}
This Fock module has a VOA structure on it, namely, one can define two quantum fields:
\begin{equation}
b(z)=\sum_m b_mz^{-m-2}, \qquad c(z)=\sum_n c_nz^{-n+1},
\end{equation}
which according to the commutation relations between modes have the following operator product:
\begin{equation}
b(z)c(w)\sim\frac{1}{z-w}.
\end{equation}
The Virasoro element is given by the following expression:
\begin{equation}
L^{\Lambda}(z)=2:\partial b(z)c(z):+:b(z)\partial c(z):,
\end{equation}
such that $b(z)$ has conformal weight $2$, and $c(z)$ has conformal weight $-1$.

The central charge of the corresponding Virasoro algebra is equal to -26.
One can define the following operator:
\begin{equation}
N_g(z)=:c(z)b(z):
\end{equation}
which is known as ghost number current.
The reason for such name is the following one. The operator $N_g=\oint \frac{\ud z}{2\pi i}N_g(z)$ gives
an integer grading to the Fock module $\Lambda$, namely,
\begin{equation}
N_g\mathbf{1}=0,\quad [N_g,b_n]=-b_n, \quad [N_g,c_m]=c_m.
\end{equation}
Next, we consider the semi-infinite cohomology (BRST) operator \cite{ko}, \cite{FGZ}. 
Let $V$ be the space, equipped with the structure of the VOA, where the Virasoro algebra has central charge,
equal to 26. 
Let's consider the tensor product $V\otimes \Lambda$. This space has a structure of VOA, such that
the central charge of the Virasoro algebra is equal to 0.
Moreover, the space $V\otimes\Lambda$ has integer grading with respect to the operator $N_g$. The following Proposition holds (see e.g. \cite{FGZ}).
\begin{prop}
The operator of ghost number 1
\begin{equation}
Q=\int \frac{\ud z}{2\pi i}J_B(z), \quad J_B(z)=:\big(L^V(z)+\frac{1}{2}L^\Lambda(z)\big)c(z):+\frac{3}{2}\partial^2c(z)
\end{equation}
is nilpotent: $Q^2=0$ on $V\otimes\Lambda$.
\end{prop}
 The space $V\otimes\Lambda$ is known as semi-infinite cohomology complex, where the differential is Q-operator
known in physics literature as BRST operator. 
The grading in the complex is given by ghost number operator $N_g$. The $k$-th cohomology group
is usually denoted as $H^{\frac{\infty}{2}+k}(Vir,\mathbb{C}{\bf{c}}, V)$.\\
The following operator product expansions will be helpful for us below (see e.g. \cite{fb}):
\begin{eqnarray}
J_B(z)b(w)&\sim&\frac{3}{(z-w)^2}+\frac{1}{(z-w)^2}N_g(w)+\nonumber\\
& &\frac{1}{z-w}\big(L^V(w)+L^\Lambda(w)\big),\nonumber\\
J_B(z)c(w)&\sim&\frac{1}{z-w}c(w)\partial c(w),\nonumber\\
J_B(z)Y(a,w)&\sim&\frac{\Delta_a}{(z-w)^2}c(w)Y(a,w)+\nonumber\\
& &\frac{1}{z-w}\big(\Delta_a\partial c(w)Y(a,w)+c(w)\partial_wY(a,w)\big).
\end{eqnarray}
As a consequence we obtain the following statement.
\begin{col}\label{qbcy}
\begin{eqnarray}
&&[Q,c(z)]=c(z)\partial c(z),\nonumber\\
&&[Q,b(z)]=L^{\Lambda}(z)+L^{V}(z),\nonumber\\
&&[Q,Y(a,z)]=\Delta_a\partial c(w)Y(a,w)+c(w)\partial_wY(a,w).
\end{eqnarray}
In the formulas above "$a$" denotes the highest weight vector w.r.t. $\{L_n^V\}$ Virasoro algebra
and $Y(a,z)$ denotes the corresponding vertex operator.
\end{col}
\noindent Finally, we remind the following crucial fact about semi-infinite cohomology.
\begin{lemma}
Let $\Phi\in V$, such that $\Phi$ is $Q$-closed and $L_0\Phi=\Delta \Phi$. Then if $\Delta\neq 0$, 
$\Phi=Q\Psi$ for some $\Psi\in V$. In other words, semi-infinite cohomology is nontrivial only on the level 
$L_0=0$.
\end{lemma}
\noindent\underline{\bf Proof.\ }
We know that $[Q,b_0]=L_0$, therefore,
\begin{eqnarray}
[Q,b_0]\Phi=Qb_0\Phi=\Delta \Phi.
\end{eqnarray}
Therefore, for $\Delta\neq 0$, $\Psi=\Delta^{-1}b_0\Phi$.
\hfill$\blacksquare$\bigskip

\noindent{\bf{6.2. Double of the braided vertex algebra and semi-infinite cohomology.}} 
We remind that above we defined the braided VOA structure on the space
$
\mathbb{F}_\varkappa=\oplus_\lambda(V_{\Delta(\lambda),\varkappa}\otimes V_{\lambda}^q)
$.
We also remind that $V_{\lambda}^q$ is the irreducible representation for $U_q(sl(2))$, where
$q=e^{\frac{\pi i}{\varkappa}}$. Let's consider the space $\mathbb{F}=\mathbb{F}_\varkappa\otimes\mathbb{F}_{-\varkappa}$.

\begin{prop}
The space $\mathbb{F}=\mathbb{F}_\varkappa\otimes\mathbb{F}_{-\varkappa}$ possesses a structure of
braided VOA such that the Virasoro algebra has central charge 26.
\end{prop}

\noindent\underline{\bf Proof.\ }
It is clear that one can define a map
\begin{equation}
\hat{Y}:\mathbb{F}_\varkappa\otimes\mathbb{F}_{-\varkappa}\rightarrow End(\mathbb{F}_\varkappa\otimes\mathbb{F}_{-\varkappa})\{z\}
\end{equation}
in the following way:
\begin{equation}
\hat{Y}\big((v\otimes a)\otimes(\bar{v}\otimes\bar{a}),z\big)=
Y(v\otimes a,z)Y(\bar{v}\otimes\bar{a},z),
\end{equation}
where $v\in F_{\lambda,\varkappa},\ a\in V_{\lambda}^q, \ \bar{v}\in F_{\mu,-\varkappa}, \ \bar{a}\in V_{\mu}^{q^{-1}}$.
It is easy to check that the map $\hat{Y}$ satisfies all properties of the braided vertex algebra, where 
the commutativity is given by:
\begin{eqnarray}
&&\mathscr{A}_{z,w}\hat{Y}\big((v_1\otimes a_1)\otimes(\bar{v}_1\otimes \bar{a}_1),z\big)
\hat{Y}\big((v_2\otimes a_2)\otimes(\bar{v}_2\otimes\bar{a}_2),w\big)=\nonumber\\
&&\hat{Y}\big((v_2\otimes r_i^{(1)}a_2)\otimes(\bar{v}_2\otimes\bar{r}_i^{(1)}\bar{a}_2),w\big)\nonumber\\
&&\hat{Y}\big((v_1\otimes r_i^{(2)}a_1)\otimes(\bar{v}_1\otimes\bar{r}_i^{(2)}a_2),z\big),
\end{eqnarray}
where $R=\sum_i r_i^{(1)}\otimes r_i^{(2)}$ is the universal R-matrix for $U_q(sl(z))$ and 
$\bar{R}=\sum_i\bar{r}_i^{(1)}\otimes\bar{r}_i^{(2)}$ is the universal R-matrix for $U_{q^{-1}}(sl(z))$.

The vertex operator, corresponding to the Virasoro element for $\mathbb{F}$, is given by
\begin{equation}
L(z)=L^\varkappa(z)+L^{-\varkappa}(z),
\end{equation}
where $L^\varkappa(z)$ and $L^{-\varkappa}(z)$ are the vertex operators associated with Virasoro elements of $\mathbb{F}_\varkappa$ and $\mathbb{F}_{-\varkappa}$ correspondingly. From Proposition \ref{feiginfuks} that the central charge, corresponding to Virasoro algebra generated by $L(z)$ is given
by
\begin{equation}
c=13-6(\varkappa+\frac{1}{\varkappa})+13+6(\varkappa+\frac{1}{\varkappa})=26.
\end{equation}
\hfill$\blacksquare$\bigskip

In such a way, $\mathbb{F}\otimes\Lambda$ also possesses a braided VOA structure. 
It is clear that  $Q$-closed terms form a braided VOA subalgebra in $\mathbb{F}\otimes\Lambda$. 
It is important to note that since the nontrivial cohomology occurs only on the level $L_0=0$, then 
nontrivial cycles with respect to $Q$ belong to the subspace $\mathbb{F}^r=\oplus_{\lambda\in \mathbb{Z}_+}\mathbb{F}_\varkappa(\lambda)\otimes\mathbb{F}_{-\varkappa}(\lambda)$, since this is the subspace containing all the zeros of $L_0$ in the generic case. Moreover, we have the following statement.
\begin{prop}\label{local}
The cycles with respect to Q in $\mathbb{F}$ involve only integer powers of formal variables in the operator products up to the cohomologically trivial terms.  
\end{prop}
\noindent\underline{\bf Proof.\ }
Let's look on the operator product
\begin{equation}
\hat{Y}\big((v_1\otimes a_1)\otimes(\bar{v}_1\otimes \bar{a}_1),z\big)
(v_2\otimes a_2)\otimes(\bar{v}_2\otimes\bar{a}_2),
\end{equation}
where $v_i\otimes a_i\in \mathbb{F}_\varkappa(\lambda_i),\  \bar{v}_1\otimes \bar{a}_1\in\mathbb{F}_{-\varkappa}(\bar{\lambda}_i)$
for some $\lambda_i,\bar{\lambda}_i\in\mathbb{Z}$.
This product belongs to the following space:
\begin{equation}
\bigoplus_{\lambda_3,\bar{\lambda}_3}\mathbb{F}_\varkappa(\lambda_3)\otimes\mathbb{F}_{-\varkappa}(\bar{\lambda}_3)
z^{\Delta_{12}(\lambda_3)}\cdot\big[[z,z^{-1}]\big],
\end{equation}
where
$\Delta_{12}(\lambda_3)=\Delta(\lambda_3)+\bar{\Delta}(\bar{\lambda}_3)-
\Delta(\lambda_1)-\Delta(\lambda_2)-\bar{\Delta}(\bar{\lambda}_1)-\bar{\Delta}(\bar{\lambda}_2)$ and
$\Delta(\lambda)=-\frac{\lambda}{2}+\frac{\lambda(\lambda+2)}{4\varkappa}$, 
$\bar{\Delta}(\bar{\lambda})=-\frac{\bar{\lambda}}{2}-\frac{\bar{\lambda}(\bar{\lambda}+2)}{4\varkappa}$.
We see that for $\varkappa\notin \mathbb{Q}$ this operator product contains a "regular'' part (i.e. just integer powers of z) only
in the case if $\lambda_1=\bar{\lambda}_1,\ \lambda_2=\bar{\lambda}_2$.
This regular part will correspond to the sector $\lambda_3=\bar{\lambda}_3$.

Hence, $\hat{Y}^r$, the reduction of $\hat{Y}$ to the regular part, is equivalent to the reduction of $\hat{Y}$ to
the subspace $\mathbb{F}^r=\oplus_\lambda\big(\mathbb{F}_\varkappa(\lambda)\otimes\mathbb{F}_{-\varkappa}(\lambda)\big)$. 
But this is precisely, what we need to do if we want to cancel the $Q$-exact terms in case if $v_1\otimes\bar{v}_1$ and $v_2\otimes \bar{v}_2$ are $Q$-closed. 
\hfill$\blacksquare$\bigskip

\noindent The last calculation we do in this section corresponds to the explicit form of the 
commutativity relation in the simplest nontrivial case, when $\lambda=1$. 
There are only two vectors in each of 
$V_1^q,\ V_1^{q^{-1}}$, i.e. the highest weight and the lowest weight vectors.
We denote them as $a_+, a_-$ and $\bar{a}_+,\bar{a}_-$ correspondingly. Let us make the following notation:
\begin{eqnarray}
v\otimes a_-\otimes\bar{v}\otimes\bar{a}_+&=&A^{(v,\bar{v})},\nonumber\\
v\otimes a_+\otimes\bar{v}\otimes\bar{a}_-&=&D^{(v,\bar{v})},\nonumber\\
v\otimes a_+\otimes\bar{v}\otimes\bar{a}_+&=&B^{(v,\bar{v})},\nonumber\\
v\otimes a_-\otimes\bar{v}\otimes\bar{a}_-&=&C^{(v,\bar{v})}
\end{eqnarray}
for any $v\in V_{\Delta(\lambda),\varkappa},\ \bar{v}\in V_{\bar{\Delta}(\lambda),-\varkappa}$.
Then the following Proposition holds.
\begin{prop}\label{abcd}
Let $v_i\otimes \bar{v}_i\in V_{\Delta(1),\varkappa}\otimes V_{\bar{\Delta}(1),-\varkappa}$ $(i=1,2)$ be $Q$-closed. The commutativity relation on $\mathbb{F}$ leads to the following  relations, which hold in 
$\mathbb{F}$ up to $Q$-exact terms:
\begin{eqnarray}\label{abcdop}
\mathscr{A}_{z,w}(B_1(z)A_2(w))&\approx&A_2(w)B_1(z)q^{-1}+(q-q^{-1})q^{-1}B_2(w)A_1(z),\nonumber\\
\mathscr{A}_{z,w}(B_1(z)C_2(w))&\approx&C_2(w)B_1(z)+
(q-q^{-1})(D_2(w)A_1(z)-\nonumber\\
&&A_2(w)D_1(z))-(q-q^{-1})^2B_2(w)C_1(z),\nonumber\\
\mathscr{A}_{z,w}(B_1(z)D_2(w))&\approx&D_2(w)B_1(z)q-B_2(w)A(z)q(q-q^{-1}) ,\nonumber\\
\mathscr{A}_{z,w}(A_1(z)C_2(w))&\approx&C_2(w)A_1(z)q-A_2(w)C_1(z)q(q-q^{-1}),\nonumber\\
\mathscr{A}_{z,w}(D_1(z)A_2(w))&\approx&A_2(w)D_1(z) +(q-q^{-1})B_2(w)C_2(z),\nonumber\\
\mathscr{A}_{z,w}(D_1(z)C_2(w))&\approx&C_2(w)D_1(z)q^{-1}+(q-q^{-1})q^{-1}D_2(w)C_1(z), \nonumber\\
\mathscr{A}_{z,w}(A_1(z)B_2(w))&\approx&B_2(w)A_1(z)q^{-1},\nonumber\\
\mathscr{A}_{z,w}(C_1(z)B_2(w))&\approx&B_2(w)C_1(z),\nonumber\\
\mathscr{A}_{z,w}(D_1(z)B_2(w))&\approx&B_2(w)D_1(z)q,\nonumber\\
\mathscr{A}_{z,w}(C_1(z)A_2(w))&\approx&A_2(w)C_1(z)q,\nonumber\\
\mathscr{A}_{z,w}(C_1(z)D_2(w))&\approx&D_2(w)C_1(z)q^{-1},\nonumber\\
\mathscr{A}_{z,w}(A_1(z)D_2(w))&\approx&D_2(w)A_1(z)+(q^{-1}-q)B_2(w)C_1(z),
\end{eqnarray}
where $S_i(S=A,B,C,D, \ i=1,2)$ stands for $S^{(v_i,\bar{v}_i)}(S=A,B,C,D, \ i=1,2)$.
\end{prop}

\noindent\underline{\bf Proof.\ }
One can obtain all these relations by means of direct use of the R-matrix. Really, in the fundamental representation
the R-matrix acts as follows:
\begin{eqnarray}
R&=&q^{\frac{H\otimes H}{2}}\big(1+(q-q^{-1})E\otimes F\big),\nonumber\\
\bar{R}&=&q^{-\frac{H\otimes H}{2}}\big(1-(q-q^{-1})E\otimes F\big),
\end{eqnarray}
since the higher powers of $E$ and $F$ act as $0$. Since we chose $a_-=Fa_+$ and $\bar{a}_-=\bar{F}\bar{a}_+$,
the result can be obtained by  direct computation. We will give here the explicit computation of the first line in (\ref{abcdop}) (all other relations can be derived in a similar way).  Namely, 
\begin{eqnarray}
&&\mathscr{A}_{z,w}(B_1(z)A_2(w))=\nonumber\\
&&\mathscr{A}_{z,w}(a_+\otimes v_1\otimes \bar{a}_+\otimes \bar{v}_1)(z)(a_-\otimes v_1\otimes \bar{a}_+\otimes \bar{v}_1)(w))\approx\nonumber\\
&&R\bar{R}(a_-\otimes v_1\otimes \bar{a}_+\otimes \bar{v}_1)(w)(a_+\otimes v_1\otimes \bar{a}_+\otimes \bar{v}_1)(z)=\nonumber\\
&&q^{-1}(a_-\otimes v_1\otimes \bar{a}_+\otimes \bar{v}_1)(w)(a_+\otimes v_1\otimes \bar{a}_+\otimes \bar{v}_1)(z)+\nonumber\\
&&(q-q^{-1})q^{-1}(a_+\otimes v_1\otimes \bar{a}_+\otimes \bar{v}_1)(w)(a_-\otimes v_1\otimes \bar{a}_+\otimes \bar{v}_1)(z)=\nonumber\\
&&A_2(w)B_1(z)q^{-1}+(q-q^{-1})q^{-1}B_2(w)A_1(z).
\end{eqnarray}
\hfill$\blacksquare$\bigskip

\noindent This Proposition will be used in the identification of the ring structure of 
 $H^{\frac{\infty}{2}+\cdot}(Vir, \mathbb{C}{\bf{c}},\mathbb{F})$ that we will now define.\\

\noindent {\bf 6.3. Lian-Zuckerman associative algebra and $SL_q(2)$.} Let's consider the following algebraic operation on the $H^{\frac{\infty}{2}+\cdot}(Vir, \mathbb{C}{\bf{c}},\mathbb{F})$:
\begin{equation}
\mu(U,V)=\mathrm{Res}_z \big(\frac{U(z)V}{z}\big).
\end{equation}
Here $U,V$  are representatives of  $H^{\frac{\infty}{2}+\cdot}(Vir, \mathbb{C}{\bf{c}},\mathbb{F})$ and $U(z)$ is shorthand notation for the vertex operator 
corresponding to the vector $U$. Since $[Q,A(z)]V=QA(z)V$ and $U(z)QB=QU(z)B$ for any $A,B\in \mathbb{F}$, one obtains that the operation $\mu$ is well defined, i.e. does not depend on the choice of representatives in $H^{\frac{\infty}{2}+\cdot}(Vir, \mathbb{C}{\bf{c}},\mathbb{F})$. 
Then the following statement holds. 
\begin{prop}
The operation $\mu$ being considered on $H^{\frac{\infty}{2}+\cdot}(Vir,\mathbb{C}\mathbf{c},\mathbb{F})$
is associative and satisfies the following commutativity relation:
\begin{equation}
\mu(U,V)=\mu(\hat{r}_i^{(1)}V,\hat{r}_i^{(2)}U)(-1)^{|U||V|},
\end{equation}
where $\hat{R}=\sum_i \hat{r}_i^{(1)}\otimes\hat{r}_i^{(2)}=R\bar{R}$ and $|\cdot|$ denotes the ghost number.
\end{prop}

\noindent\underline{\bf Proof.\ }
 The proof follows the same steps as in "abelian" case (nonbraided VOA) studied in \cite{lz2}. Due to the 
Proposition \ref{local} we can limit ourselves to the consideration only of the regular terms in operator products (all other contributions correspond to the $Q$-exact terms). At first, let us prove commutativity relation:
\begin{eqnarray}
&&\mu(U,V)-(-1)^{|U||V|}\mu(\hat{r}_i^{(1)}V,\hat{r}_i^{(2)}U)=\nonumber\\
&&\mathrm{Res}_w\mathrm{Res}_{z-w}\frac{\big(U(z-w)V\big)(w)\mathbf{1}}
   {\big(1+(z-w)/w\big)w^2}=\nonumber\\
&&\sum_{i\geqslant 0}(-1)^i\mathrm{Res}_w\mathrm{Res}_{z-w}
  \frac{\big(U(z-w)V\big)(w)\mathbf{1}}{(z-w)^{-i}w^{i+2}}=\nonumber\\
&&\sum_{i\geqslant 0}\frac{(-1)^i}{i+1}\mathrm{Res}_w\mathrm{Res}_{z-w}L_{-1}
   \frac{\big(U(z-w)V\big)(w)\mathbf{1}}{(z-w)^{-i}w^{i+1}}=\nonumber\\
&&\sum_{i>0}\frac{(-1)^i}{i+1}\mathrm{Res}_w\mathrm{Res}_{z-w}
   \frac{(Qb_{-1}+b_{-1}Q)\big(U(z-w)V\big)(w)\mathbf{1}}{(z-w)^{-i}w^{i+1}}=\nonumber\\
&&\sum_{i>0}Q\mathrm{Res}_w\mathrm{Res}_{z-w}b_{-1}
   \frac{\big(U(z-w)V\big)(w)\mathbf{1}}{(z-w)^{-i}w^{i+1}}
\end{eqnarray}
for any $U,V$ which are the representatives of the cohomology classes 
$H^{\frac{\infty}{2}+\cdot}(Vir,\mathbb{C}\mathbf{c},\mathbb{F})$. Therefore, the commutativity relation
holds. Now, let us prove the associativity of $\mu$.
\begin{eqnarray}
&&\mu\big(\mu(U,V),W\big)-\mu\big(U,\mu(V,W)\big)=\nonumber\\
&&\mathrm{Res}_w\mathrm{Res}_{z-w}
     \frac{\big(U(z-w)V\big)(w)W}{(z-w)w}-\nonumber\\
&&   \mathrm{Res}_z\mathrm{Res}_w\frac{U(z)V(w)W}{zw}=
  \mathrm{Res}_z\mathrm{Res}_w\frac{U(z)V(w)W}{(1-\frac{w}{z})zw}+\nonumber\\
&&   (-1)^{|U||V|}\sum_k\mathrm{Res}_w\mathrm{Res}_z
   \frac{(\hat{r}_k^{(1)}V)(w)(\hat{r}_k^{(2)}U)(z)W}{(1-\frac{z}{w})w^2}-\nonumber\\
&&   \mathrm{Res}_z\mathrm{Res}_w\frac{U(z)V(w)W}{zw}=
\sum_{i>0}\mathrm{Res}_z\mathrm{Res}_w\frac{U(z)V(w)W}{z^{i+1}w^{-i+1}}+\nonumber\\
&&(-1)^{|U||V|}\sum_k\sum_{i\geqslant 0}\mathrm{Res}_w\mathrm{Res}_z
     \frac{(\hat{r}_k^{(1)}V)(w)(\hat{r}_k^{(2)}U)(z)W}{z^{-i}w^{i+2}}=\nonumber\\
&&\sum_{i\geqslant 0}\frac{1}{i+1}\mathrm{Res}_z\mathrm{Res}_w
     \frac{(L_{-1}U)(z)V(w)W}{z^{i+1}w^{-i}}+\nonumber\\
&&(-1)^{|U||V|}\sum_k\sum_{i>0}\frac{1}{i+1}\mathrm{Res}_w\mathrm{Res}_z
     \frac{(L_{-1}\hat{r}_k^{(1)}U)(w)\hat{r}_k^{(2)}V(z)W}{z^{-i}w^{i+1}}=\nonumber\\
&&Qn(U,V,W)+n(QU,V,W)+(-1)^{|U|}n(U,QV,W)+\nonumber\\
&&(-1)^{|U|+|V|}n(U,V,QW),
\end{eqnarray}
where $n$ is given by
\begin{eqnarray}
&&n(U,V,W)=\sum_{i\geqslant 0}\frac{1}{i+1}\mathrm{Res}_z\mathrm{Res}_w
     \frac{(b_{-1}U)(z)V(w)W}{z^{i+1}w^{-i}}+\nonumber\\
&&(-1)^{|U|+|V|}\sum_{i\geqslant 0}\frac{1}{i+1}\mathrm{Res}_w\mathrm{Res}_z
     \frac{(b_{-1}\hat{r}_k^{(1))}U)(z)(\hat{r}_k^{(2)}V)(w)W}{z^{-i}w^{i+1}}.
\end{eqnarray}
Therefore, $\mu$ is associative on semi-infinite cohomology. Thus, the Proposition is proven.
\hfill$\blacksquare$\bigskip

The next Proposition gives a possibility to compute the semi-infinite cohomology of $\mathbb{F}$.

\begin{prop}
Let $\mathbf{F}=\oplus_{\lambda\in \mathbb{Z}_+}(V_{\Delta(\lambda),\varkappa}\otimes V_{\bar{\Delta}(\lambda),-\varkappa})$.
Then \\
i) $\mathbf{F}$ has a VOA structure;\\
ii) $H^{\frac{\infty}{2}+k}\big(Vir,\mathbb{C}\mathbf{c},\mathbf{F}(\lambda)\big)=\mathbb{C}\delta_{k,0}\oplus\mathbb{C}\delta_{k,3}$,
where $\mathbf{F}(\lambda)=V_{\Delta(\lambda),\varkappa}\otimes V_{\bar{\Delta}(\lambda),-\varkappa}$.
\end{prop}
\noindent The proof is given in \cite{frst}, following the results of \cite{lz1}, \cite{lz1a}. It is not hard to calculate the explicit formulas for the representatives of 
$H^{\frac{\infty}{2}+0}\big(Vir,\mathbb{C}\mathbf{c},\mathbf{F})$ and 
$H^{\frac{\infty}{2}+3}\big(Vir,\mathbb{C}\mathbf{c},\mathbf{F})$.

\begin{prop}
\hfill
\begin{itemize}
\item[(i)]
The operators corresponding to the representatives of $H^{\frac{\infty}{2}+0}\big(Vir,\mathbb{C}\mathbf{c},\mathbf{F}(1))$ 
and $H^{\frac{\infty}{2}+3}\big(Vir,\mathbb{C}\mathbf{c},\mathbf{F}(1))$ have the following explicit form:
\begin{eqnarray}
\Phi^0(z)&=&L^{\varkappa}_{-1}\Phi(z)-L^{-\varkappa}_{-1}\Phi(z)-\varkappa^{-1}:bc:(z)\Phi(z),\nonumber\\
\Phi^3(z)&=&c\partial c\partial^2 c L^{tot}_{-1}\Phi(z),
\end{eqnarray}
where we denoted by $L^{\varkappa}_n,L^{-\varkappa}_n$ the Virasoro algebra generators in 
$V_{\Delta(\lambda),\varkappa}$, \\*
$V_{\bar{\Delta}(\lambda),-\varkappa}$ correspondingly and
$L_n\equiv L^{\varkappa}_n+L^{-\varkappa}_n$.
The operator $\Phi(z)$ corresponds to the tensor product of highest weight states in
$V_{\Delta(1),\varkappa}$ and $V_{\bar{\Delta}(1),-\varkappa}$;
\item[(ii)]$\Phi^0,\mathbf{1}$ (the vector corresponding to $\Phi^0(z)$ and the vacuum state) generate
all $H^{\frac{\infty}{2}+0}\big(Vir,\mathbb{C}\mathbf{c},\mathbf{F})$ by means of bilinear operation $\mu$.
\end{itemize}
\end{prop}
\noindent\underline{\bf Proof.\ }
\\*
{\bf (i)} First of all, we have the following relation, which follows from (\ref{qbcy}):
\begin{equation}
[Q,\Phi](z)=-\partial c\Phi+c\partial\Phi.
\end{equation}Let $L^{\varkappa}(z)$ be the Virasoro element in $V_{\Delta(\lambda),\varkappa}$. Then 
$L^{\varkappa}_{-1}\Phi(w)=\int\frac{\ud z}{2\pi i}L^{\varkappa}(z)\Phi(w)$. Therefore,
\begin{equation}
[Q,L^{\varkappa}_{-1}\Phi](w)=\int\frac{\ud z}{2\pi i}L^{\varkappa}(z)[Q,\Phi](w)=\int\frac{\ud z}{2\pi i}[Q,L^{\varkappa}(z)]\Phi(w).
\end{equation}
It is easy to see that 
\begin{equation}
[Q,\int\frac{\ud z}{2\pi i}L^{\varkappa}(z)]=\int\frac{\ud z}{2\pi i}2\partial^2cL^{\varkappa}(z)+c\partial L^{\varkappa}(z).
\end{equation}
Hence,
\begin{eqnarray}
\int\frac{\ud z}{2\pi i}[Q,L^{\varkappa}(z)]\Phi(w)&=&\Delta(1)\partial^2 c\Phi+\partial cL_{-1}^{\varkappa}\Phi,\nonumber\\
\int\frac{\ud z}{2\pi i}L^{\varkappa}(z)[Q,\Phi](w)&=&-\partial c\bar{L}^{\varkappa}_{-1}\Phi+cL^{-\varkappa}_{-1}\partial\Phi,
\end{eqnarray}
which lead to
\begin{equation}
[Q,L^{\varkappa}_{-1}\Phi(w)]=\Delta(1)\partial^2 c\Phi+c{L}^{-\varkappa}_{-1}\partial\Phi.
\end{equation}
It is not hard to see that
\begin{equation}
[Q,(L^{\varkappa}_{-1}-L^{-\varkappa}_{-1})\Phi(w)]=
\frac{3}{2\varkappa}\partial^2 c\Phi+c(({L_{-1}^{-\varkappa}})^2-(L_{-1}^{\varkappa})^2)\Phi.
\end{equation}
Another necessary fact which we will need is that there is a linear dependence between 
$({L_{-1}^{\varkappa}})^2\Phi$ and $L_{-2}^{-\varkappa}\Phi$, namely,
\begin{eqnarray}
({L^{\varkappa}_{-1}})^2\Phi&=&\frac{2}{3}\big(1+2\Delta(1)\big)L^{\varkappa}_{-2}\Phi,\nonumber\\
({L_{-1}^{-\varkappa}})^2\Phi&=&\frac{2}{3}\big(1+2\bar{\Delta}(1)\big)L^{-\varkappa}_{-2}\Phi
\end{eqnarray}
or
\begin{equation}
({L^{\varkappa}_{-1}})^2\Phi=\varkappa^{-1}L^{\varkappa}_{-2}\Phi,\qquad ({L_{-1}^{-\varkappa}})^2\Phi=-\varkappa^{-1}L^{-\varkappa}_{-2}\Phi.
\end{equation}
Therefore,
\begin{equation}
[Q,(L^{\varkappa}_{-1}-L_{-1}^{-\varkappa})\Phi(w)]=\varkappa^{-1}\big(\frac{3}{2}\partial^2 c\Phi+
cL_{-2}\Phi\big).
\end{equation}
Let us consider another term from $\Phi^0$, namely, $:bc:(w)\Phi(w)$
\begin{eqnarray}
[Q,:bc:(w)\Phi(w)]&=&\int_w\frac{\ud z}{2\pi i}\frac{1}{(z-w)}\big([Q,:bc:](z)\Phi(w)\nonumber\\
&+&:bc:(z)\int [Q,\Phi(w)]\big),\nonumber\\
\lbrack Q,:bc:(w)\rbrack=j_B(w)&=&c(w)L(w)+:c\partial cb:(w)+\frac{3}{2}\partial^2 c(w),\nonumber\\
\int_w\frac{\ud z}{2\pi i}\frac{j_B(z)\Phi(w)}{z-w}&=&-\frac{1}{2}\partial^2 c\Phi+
\partial c\partial\Phi+cL_{-2}\Phi\nonumber\\
&+&:bc\partial c:\Phi+\frac{3}{2}\partial^2 c\Phi.
\end{eqnarray}
On the other hand,
\begin{eqnarray}
\int_w\frac{\ud z}{2\pi i}\ \frac{:bc:(z)}{z-w}\big(-\partial c(w)\Phi(w)+c(w)\partial\Phi(w)\big)=\nonumber\\
\frac{1}{2}\partial^2 c(w)\Phi(w)-\partial c(w)\partial\Phi(w)-:bc\partial c:(w)\Phi(w).
\end{eqnarray}
Therefore,
\begin{equation}
[Q,:bc:(w)\Phi(w)]=\frac{3}{2}\partial^2 c\Phi(w)+cL_{-2}\Phi(w),
\end{equation}
and the operator
\begin{equation}
(L^{\varkappa}_{-1}-{L_{-1}^{-\varkappa}})\Phi(w)-\varkappa^{-1}:bc:(w)\Phi(w)
\end{equation}
corresponds to the closed state in ${\bf F}(1)\otimes \Lambda$. At the same time, it is the only $Q$-closed state of ghost number $0$ in 
${\bf F}(1)\otimes \Lambda$, which belongs to the kernel of $L_0$, therefore, due to Theorem 3.7 it is not exact. 
Hence, it represents the cohomology class from $H^{\frac{\infty}{2}+0}(Vir,\mathbb{C}\mathbf{c},\mathbf{F})$.
It is not hard to check that $\Phi^3=c\partial c\partial^2 cL_{-1}\Phi(z)$
represents the cohomology class $H^{\frac{\infty}{2}+3}(Vir,\mathbb{C}\mathbf{c},\mathbf{F})$.
The part \noindent{\bf (ii)} follows from the Theorem 3.7. of \cite{frst}. 
\hfill$\blacksquare$\bigskip

\noindent Our next aim is to use the above information to calculate 
$H^{\frac{\infty}{2}+\cdot}(Vir,\mathbb{C}\mathbf{c},\mathbb{F})$ and 
to understand its underlying multiplicative structure.

\begin{prop}
The 0-th and 3-d semi-infinite cohomology groups of $\mathbb{F}$ are
\begin{eqnarray}
H^{\frac{\infty}{2}+0}(Vir,\mathbb{C}\mathbf{c},\mathbb{F})\cong
H^{\frac{\infty}{2}+3}(Vir,\mathbb{C}\mathbf{c},\mathbb{F})\cong
\bigoplus_{\lambda\geqslant 0}V^q_\lambda\otimes V_\lambda^{q^{-1}}.
\end{eqnarray}
\end{prop}

\noindent\underline{\bf Proof.\ }
The proof evidently follows from the fact:
\begin{equation}
H^{\frac{\infty}{2}+0}(Vir,\mathbb{C}\mathbf{c},\mathbf{F}(\lambda))\cong
H^{\frac{\infty}{2}+3}(Vir,\mathbb{C}\mathbf{c},\mathbf{F}(\lambda))\cong\mathbb{C}
\end{equation}
for $\lambda\ge 0$.
\hfill$\blacksquare$\bigskip

\begin{theorem}
$H^{\frac{\infty}{2}+0}\big(Vir,\mathbb{C}\mathbf{c},\mathbb{F}(1)\big)$ generates
all $H^{\frac{\infty}{2}+0}(Vir,\mathbb{C}\mathbf{c},\mathbb{F})$ by means of multiplication
$\mu$, and the generating set 
$${A,B,C,D}\in H^{\frac{\infty}{2}+0}\big(Vir,\mathbb{C}\mathbf{c},\mathbb{F}(1)\big)
\cong V_1^q\otimes V_1^{q^{-1}}
$$ 
satisfies the following relations:
\begin{eqnarray}
&&AB=BA{q^{-1}}, \quad CB=BC,\quad DB=BDq,\quad CA=ACq,\nonumber\\
&&AD-DA=(q^{-1}-q)BC,\quad CD=DCq^{-1},\nonumber\\
&&AD-q^{-1}BC=1,
\end{eqnarray}
or, equivalently,
\begin{equation}
\big(H^{\frac{\infty}{2}+0}(Vir,\mathbb{C}\mathbf{c},\mathbf{F}),\mu\big)\cong SL_q(2).
\end{equation}
\end{theorem}

\noindent\underline{\bf Proof.\ }
The statement, that 
$H^{\frac{\infty}{2}+0}(Vir,\mathbb{C}\mathbf{c},\mathbb{F}(1))=V_1^q\otimes V_1^{q^{-1}}$
forms a generating set for 
$H^{\frac{\infty}{2}+0}(Vir,\mathbb{C}\mathbf{c},\mathbb{F})=\oplus_{\lambda\in \mathbb{Z}_{\ge 0}} V_\lambda^q\otimes V_\lambda^{q^{-1}}$, 
is a consequence of the fact that by means of tensor multiplication of $V_1^q$ with each other one can obtain any of $V_\lambda^q$ and
the explicit form of vertex operators (since they act as intertwiners).
One can choose generators $A,B,C,D$ in the same way as we did in
Proposition \ref{abcd}. Then all the relations except for the last one follow from the commutativity relation of $\mu$. One can obtain the last relation by means of 
the explicit consideration of the action of vertex operators (more precisely, corresponding intertwiners) and normalizing $A,B,C,D$ appropriately. 
\hfill$\blacksquare$\bigskip

\noindent In this section, we found how the multiplication structure on $SL_q(2)$ emerges from braided VOA structure on $\mathbb{F}$ via Lian-Zuckerman construction.
However, $SL_q(2)$ is a Hopf algebra and the natural question is whether comultiplication structure has the same origin. The answer is  yes, and it is provided by means of the notion of vertex operator coalgebra (VOCA) \cite{hub}. 
If one has a map $Y$ giving a structure of VOA algebra on some space $V$ with a nondegenerate bilinear pairing $($ , $)$,  preserving the Virasoro action, the structure of VOCA is given by $Y^c:V\to V\otimes V[[z,z^{-1}]]$, such that $(Y^c(z)u,v\otimes w)=(u,Y(v,z)w)$, where $u,v,w\in V$  \cite{hub}. One can extend this definition to the case of braided VOA, namely, define the structure of a braided  VOCA algebra by the same formula. Therefore, the definition of the Lian-Zuckerman operation $\mu$ gives rise to the dual object $\Delta_{\mu}$:
\begin{eqnarray}  
&&(\Delta_{\mu}(u), v\otimes w)\equiv Res_z z^{-1}(Y^c(z)u,v\otimes w)=\nonumber\\
&&Res_z z^{-1}((u,Y(v,z)w)=(u,\mu(v,w)).
 \end{eqnarray}  
 Since $\Delta_{\mu}$ is dual to $\mu$ with respect to the canonical pairing, we obtain that $\Delta_{\mu}$ in the specific case of braided VOA $\mathbb{F}$ gives the comultiplication structure on $SL_q(2)$.  We leave the proof of the consistency of multiplication and comultiplication as an exercise. 

\section{Further developments and conjectures}

We already discussed in the introduction that thanks to the equivalence of categories of Virasoro and 
$\widehat{sl}(2)$ Lie algebras one can replace the Virasoro modules by their $\widehat{sl}(2)$
counterparts and obtain the quantum group as semi-infinite cohomology of  $\widehat{sl}(2)$ Lie
algebra. Also the extension to arbitrary affine and $\mathcal{W}$-algebras is possible thanks to 
the general results of Varchenko \cite{varchenko}, \cite{varchenko2} and Styrkas \cite{dissert}, though some 
technical difficulties might arise, related e.g. to the definition of semi-infinite cohomology of 
$\mathcal{W}$-algebras, etc (see e.g. \cite{bow}). However, we believe that the most interesting further developments will appear not 
from generalizations but from the deepening of our construction, and better understanding its 
relation with physics.

The first possible extension of our construction comes from the existence of double lattice 
(\ref{virsing}) of Verma modules with singular 
vectors and fixed central charge. The semi-infinite cohomology of the Virasoro algebra with 
values in the corresponding double lattice of the Fock spaces must yield the modular double
of the quantum group introduced in \cite{Fadd}.

The second development of our construction results from the special features 
of rational level case that we consciencely avoided in the present paper. It is 
well known that the quantum group at roots of unity has a remarkable
finite-dimensional subquotient. Its realization via semi-infinite cohomology might require
finding a corresponding subquotient braided VOA. The analogous modified regular VOA
should give rise to the finite-dimensional Verlinde algebra as conjectured in \cite{frst}.

The third possible development follows from the replacements of simple modules of 
Virasoro algebra and quantum group by the big projective modules as in \cite{frst}.
In fact, the braided VOA $\mathbb{F}_\varkappa$ has the appropriate size for such an extension 
but has a degenerate structure resulting from the absence of linking of $\mathbb{F}_{\varkappa}(\lambda)$
and $\mathbb{F}_{\varkappa}(-\lambda-2)$. The way to achieve such linking is suggested by the construction of 
the heterogeneous VOA in \cite{frst}, and will require an additional exponential term in
the action of the Virasoro algebra. 
This exponential term and the free bosonic realization of the Virasoro algebra, known as
Coulomb gas in the physics literature, reveal the relation of our construction to 
the quantum Liouville model and minimal string theory. In particular, the latter
theory also uses the pairing of representations of Virasoro algebra with 
complementary central charges $c+\bar{c}=26$.
We conjecture, that our braided VOA and its projective counterpart discussed above, admits a geometric
(and path integral) interpretation, so that the semi-infinite Lie algebra cohomology acquires
 a natural meaning as cohomology in certain infinite dimensional geometry as well. 
Then the quantum group will finally find its truly geometric setting.

\section*{Acknowledgements}
We are indebted to P.I. Etingof, Y.-Z. Huang, A.A. Kirillov Jr., M.M. Kapranov, J. Lepowsky and  G.J. Zuckerman for fruitful discussions. I.B.F. thanks K. Styrkas for a joint research on braided vertex operator algebras 10 years ago. 
A.M.Z. also wants to thank the organizers of the Simons Workshop 2008,  where this work was partly done. The research of I.B.F. was supported by NSF grant DMS-0457444.

\end{document}